\documentclass[hidelinks,onefignum,onetabnum,onealgnum]{siamart251216}

\usepackage{amsfonts,amssymb}

\usepackage{booktabs,multirow,tabularx,ragged2e}
\usepackage{graphicx,float,placeins}
\usepackage{needspace}
\usepackage{tikz}
\usetikzlibrary{arrows.meta,calc,decorations.pathreplacing}
\usepackage[algo2e,linesnumbered,ruled,vlined]{algorithm2e}
\usepackage{algpseudocode}
\newsiamremark{remark}{Remark}
\newsiamremark{example}{Example}
\newsiamremark{assumption}{Assumption}
\headers{RAS--PML for Helmholtz}{Y. Xie et al.}
\usepackage[normalem]{ulem}

\newcommand{\B}{B}
\newcommand{\Bt}{\widetilde B}

\title{Massively parallel Schwarz methods for the high-frequency Helmholtz equation with $\mathcal{O}(\lowercase{k})$ convergence}
\author{Yan Xie\thanks{Shenzhen International Center for Industrial and Applied Mathematics,
Shenzhen Research Institute of Big Data, Shenzhen, Guangdong 518172, China.}
\and Yue Yu\protect\footnotemark[1]\kern0.4em\thanks{School of Science and Engineering, The Chinese University of Hong Kong,
Shenzhen, Guangdong 518172, China.}
\and Shihua Gong\protect\footnotemark[1]\kern0.4em\protect\footnotemark[2]\kern0.4em\thanks{Shenzhen Loop Area Institute,
Shenzhen, Guangdong 518172, China. Corresponding author
(\email{gongshihua@cuhk.edu.cn}).}
\and Ivan G. Graham\thanks{Department of Mathematical Sciences, University of Bath,
Bath BA2 7AY, United Kingdom.}
\and Chen-song Zhang\thanks{SKLMS, Academy of Mathematics and Systems Science,
Chinese Academy of Sciences, Beijing 100190, China; School of Mathematical Sciences,
University of Chinese Academy of Sciences, Beijing 100049, China.}}
\begin{document}
\maketitle

\begin{abstract}
We develop a one-level restricted additive Schwarz solver with
perfectly matched layer transmission conditions (RAS--PML) for
high-frequency Helmholtz systems. For problems with $\mathcal{O}(k^d)$
unknowns, $d=2,3$, the solver uses Cartesian decompositions {and}
$\mathcal{O}(k^d)$ concurrent local solves for a single right-hand side,
where $k$ is the wavenumber. Starting from the convergence theory for fixed
decompositions established in~\texttt{arXiv:2404.02156}, we address the
growing cost {as $k$ increases} of local solves through a nested RAS--PML construction. Under
suitable stability and contraction assumptions, we prove that
sufficiently accurate recursive local solves preserve convergence of the
outer iteration. We also prove that a single inner iteration at each
non-terminal level yields a correction algebraically equivalent to
one-level RAS--PML, provided that the hierarchy-induced finest local
problems and composite restriction and weighted extension operators are
retained. This equivalence provides a hierarchical interpretation and
permits a direct implementation without recursive transfers. We assess
the one-level solver using PML and overlap widths of order
$\mathcal{O}(k^{-1}\log k)$. Constant-coefficient experiments cover
$300\leq k\leq19200$ in two dimensions and $60\leq k\leq780$ in three
dimensions, reaching $1.47\times10^9$ and $4.75\times10^8$ complex unknowns,
respectively. Richardson iteration counts grow approximately linearly
with $k$ over these ranges, while total parallel runtime, including setup
and solution, exhibits near-linear growth at the larger tested
wavenumbers. Further experiments assess the effects of mesh refinement
and heterogeneous media.
\end{abstract}



\section{Introduction}
{In one-level Schwarz methods, information is exchanged only between
neighbouring subdomains in each iteration, so the number of iterations
is tied to the number of subdomains. With the exact Dirichlet-to-Neumann
map as the subdomain boundary condition, the parallel overlapping Schwarz
method on a strip decomposition into $N$ subdomains converges in $N$
iterations~\cite{nataf1994optimal,engquist1998absorbing}. For Helmholtz problems with PML
subdomain boundary conditions, \cite{galkowski2024convergence} proves
that the number of iterations is governed by the maximum number of
subdomains between which a geometric-optic ray can travel.

For positive-definite problems, adding a coarse solve provides a global
transfer of information and makes the iteration counts independent of the
number of subdomains~\cite{toselli2004domain,dolean2015introduction}.
Recent advances in two-level Schwarz methods have also established
rigorous convergence results for high-frequency Helmholtz
problems~\cite{hu2024novel,lu2025two,ma2025two,galkowski2026convergence,graham2025two,galkowski2026twolevel,fu2026discrete}.
These analyses prove that, under their respective assumptions,
two-level methods with subdomains of size $k^{-1}$ converge in a number
of iterations independent of $k$ (possibly up to factors of $\log k$),
provided that the coarse-space dimension is at least of order $k^d$
(again, possibly up to logarithmic factors).
For example, such results are proved for problem-adapted multiscale
coarse spaces in~\cite{lu2025two,ma2025two} and $hp$-FEM coarse spaces
in~\cite{graham2025two}.
Furthermore, the numerical experiments in~\cite{galkowski2026twolevel}
show that, at least for piecewise-polynomial coarse spaces,
a pollution-free coarse FEM solution (in the sense that its relative
error is uniformly bounded in $k$) is \emph{necessary} for the number
of iterations to remain uniformly bounded in $k$.
Thus, even when the outer iteration count is bounded independently of $k$,
applying the coarse correction requires solving a Helmholtz system with at least 
$\mathcal{O}(k^d)$ unknowns.

Therefore, one cannot get away from the issue of solving the one-level
problem efficiently: the coarse problem in a two-level method is itself a
large Helmholtz system, and the same $\mathcal{O}(k^d)$ problem size}
arises in stand-alone Helmholtz discretizations with a fixed number of
degrees of freedom per wavelength in each coordinate direction.
{For a single right-hand side, efficient parallel solution requires
sufficiently many subproblems that can be solved concurrently,
with manageable local solution costs.}
We therefore focus on a one-level
solver that exposes $\mathcal{O}(k^d)$ parallelism for Helmholtz systems with
$\mathcal{O}(k^d)$ unknowns, with potential use both as a stand-alone method
and as a coarse solver within two-level methods.

\subsection{Related Work}
Domain decomposition methods {for Helmholtz problems} combine local solves with mechanisms for
transmitting wave information between subdomains. The sweeping
preconditioners of Engquist and
Ying~\cite{engquist2011sweeping,engquist2011sweeping2} approximate layerwise
factorizations, while the source-transfer method of Chen and
Xiang~\cite{chen2013source} transfers sources successively through a layered
decomposition. The dependencies between successive layers constrain
concurrent subdomain solves for a single right-hand side. Leng and
Ju~\cite{ju2021diagonal,leng2022trace} organize source or trace transfer into
diagonal sweeps on Cartesian block decompositions. Subdomains on the same
diagonal wavefront can be solved concurrently, while successive wavefronts
retain directional dependencies. For multiple right-hand sides, pipelining
can further improve processor utilization by interleaving their sweeps.

Taus et al.~\cite{taus2020sweeps} developed the L-sweeps preconditioner to
increase concurrency for a single right-hand side. At fixed resolution per
wavelength, its row-based implementation uses $\mathcal{O}(k)$ processes.
Combining the paper's cost model with the empirically observed
$\mathcal{O}(\log k)$ Krylov iteration growth gives parallel solution costs
of $\mathcal{O}(k\log k)$ in two dimensions and
$\mathcal{O}(k^2\log k)$ in three dimensions. 
The latter estimate uses
beam-shaped local problems solved sequentially on individual processes;
the authors also discuss parallelizing these local solves to reduce the
three-dimensional cost. 

Leng and Ju~\cite{leng2019additive} developed an additive overlapping
method in which all local PML problems are solved concurrently and
directional residuals are transferred between neighbouring subdomains.
Their work provides a two-dimensional constant-medium analysis and
numerical experiments in both two and three dimensions. When the number
of subdomains in each coordinate direction grows proportionally to $k$,
this formulation {allows
$\mathcal{O}(k^d)$ local problems to be solved concurrently.} The construction
of the directional residuals for source transfer requires careful design.

{Iteration counts and the cost per iteration are also central to
practical performance. Cotter et al.~\cite{cotter2025shifted} use $\mathcal{O}(k)$ iterations of a
Hermitian/skew-Hermitian splitting (HSS) method, with the resulting
complex shifted systems solved using multigrid.} Martin et al.~\cite{martin2026comparison} compare optimized restricted additive
Schwarz and substructured non-overlapping methods on large-scale
three-dimensional Helmholtz problems with multiple sources. With tuned
transmission conditions, their non-overlapping methods achieve lower
runtimes despite requiring more iterations, illustrating the importance of local factorization, solve, and Krylov costs alongside iteration counts.

\subsection{{The contributions of the present paper}}

We develop a one-level RAS--PML solver for high-frequency
Helmholtz systems with $\mathcal{O}(k^d)$ unknowns, $d=2,3$.
Our computational target is to combine $\mathcal{O}(k^d)$ parallelism
with $\mathcal{O}(k)$ parallel solution time for a single right-hand side.

{{In this work}, our starting point is RAS-PML. For a fixed domain decomposition with
$k$-independent overlap and PML widths, previous
analyses~\cite{galkowski2024convergence,galkowski2024schwarz} establish,
under nontrapping assumptions, superalgebraic decay {with respect to} $k$ of the
continuous iteration error after a fixed {(but geometry-dependent)} number of
iterations. {Thus the solver is increasingly accurate as $k$ increases. However, this} favorable high-frequency convergence must be accompanied
by an efficient treatment of the local problems. With a fixed number of
degrees of freedom per wavelength, each subdomain of fixed physical size
contains $\mathcal{O}(k^d)$ unknowns, and controlling pollution can require
further refinement \cite{galkowski2024hpfem,bernkopf2025wavenumber,li2019fem}. The computational and memory costs of direct local
solves therefore grow with frequency, making these solves a central
obstacle to a practical implementation.}

{
We address this difficulty through a nested RAS–PML construction that recursively decomposes and solves large local problems. With one RAS–PML iteration per non-terminal local solve, the correction coincides with that of one-level RAS–PML on the finest partition. This equivalence holds for the hierarchy-induced local problems and transfer operators.
It provides a hierarchical interpretation of {the method as a 
one-level residual-correction} that can be implemented without
recursive transfers. We assess the convergence and parallel performance
of the one-level solver numerically.}
{The main contributions are:}
\begin{enumerate}
\item \textbf{One-level solver and parallel performance.}
We implement a residual-correction iteration with
$\mathcal{O}(k^d)$ concurrent local PML solves, driven by restrictions of
{the residual of the global linear system}. {Experiments with growing frequency} in two and three
dimensions show approximately linear growth of the Richardson iteration
counts with $k$ and quantify the setup and solution costs of the parallel
implementation{, with wall clock time for 2D and 3D experiments growing roughly linearly in
$k$ as well.}
{For the highest frequency tested, the linear systems have}
$1.47\times10^9$ complex unknowns in two dimensions and
$4.75\times10^8$ in three dimensions. Width-selection, mesh-refinement,
and heterogeneous-media experiments further assess the parameter choices
and applicability of the method.
\item \textbf{Theory-driven recursive local solution.}
We introduce a nested RAS--PML construction in which non-terminal local
problems are themselves solved by RAS--PML iterations. 
{We assume suitable stability bounds and uniform contraction
over a fixed number of steps for the iterations with exact child solves.
We prove that sufficiently accurate recursive local solves preserve
this contraction for the nested method, with a hierarchy depth of
$\mathcal{O}(\log k)$.}
\item \textbf{Algebraic equivalence and direct implementation.}
We prove that Nested(1), with one RAS--PML iteration at each non-terminal
local level, is algebraically equivalent to a one-level method on the
finest partition when the hierarchy-induced local problems, restriction
operators, and weighted extensions are retained. 
\end{enumerate}

\subsection{Paper Organization}
The remainder of the paper is organized as follows.
Section~\ref{sec:ras} formulates the Helmholtz problem and reviews
the RAS--PML method.
{Section~\ref{sec:nested} introduces the nested RAS--PML algorithm,
states its algebraic equivalence and conditional convergence results,
and then develops their proofs.}
Section~\ref{sec:experiments}
presents numerical experiments on the convergence and parallel
performance of the one-level solver. Finally,
Section~\ref{sec:conclusion} concludes the paper.

\section{Problem formulation and the RAS--PML method}
\label{sec:ras}

\subsection{The Helmholtz problem}
Let $f\in L^2_{\mathrm{comp}}(\mathbb{R}^d)$ be a source and let
$c\in C^\infty(\mathbb{R}^d)$ be a strictly positive normalized wave
speed, with $\operatorname{supp}(1-c)$ compact. For a real reference
wavenumber $k\gg1$, we consider
{
\begin{equation}
    \label{eq:helmholtz0}
    k^{-2}\Delta u + c(x)^{-2}u = -f
    \quad\text{in}\quad\mathbb{R}^d,
\end{equation}
}
subject to the Sommerfeld radiation condition
$$\frac{\partial u}{\partial r} - \mathrm{i}k u
    = o\left(r^{-(d-1)/2}\right),
    \qquad r=|x|\to\infty.$$

We assume that the source and the medium inhomogeneity are contained in a
fixed hyperrectangle $
  \operatorname{supp}f\cup\operatorname{supp}(1-c)
  \subset \Omega_{\mathrm{int}}
  :=\prod_{i=1}^{d}(a_i,b_i).$
To truncate the unbounded problem, we surround $\Omega_{\mathrm{int}}$ by a {global}
Cartesian perfectly matched layer (PML) of thickness $\kappa_{\rm g}$ and set
$
  \Omega:=\prod_{i=1}^d(a_i-\kappa_{\rm g},b_i+\kappa_{\rm g}).
$ 
We introduce a smooth real-valued scaling function
$g\in C^{\infty}(\mathbb{R}^d)$ satisfying
{$g(x)=g'(x)=0$ for $x\leq0$ and $g'(x)>0$ for $x>0$.}
For each coordinate direction $i$, define
{
\begin{equation*}
    g_i(x_i) = \begin{cases}
        g(x_i-b_i), & x_i \geq b_i, \\
        0, & x_i \in (a_i, b_i), \\
        -g(a_i-x_i), & x_i \leq a_i.
    \end{cases}
\end{equation*}
}
The corresponding complex-scaled Laplacian is
{
\begin{equation*}
    \Delta_s := \sum_{i=1}^d
    \left(\gamma_i(x_i)^{-1}\partial_{x_i}\right)^2,
    \qquad \gamma_i:=1+\mathrm{i}g'_i.
\end{equation*}
}
The truncated PML problem is: find $u\in H^1_0(\Omega)$ such that
\begin{equation}
    \label{eq:weak}
    \begin{aligned}
    a(u,v)&:=
    \int_{\Omega}\Big(
    k^{-2}(D\nabla u)\cdot\nabla\overline v
    -k^{-2}(\beta\cdot\nabla u)\overline v
    -c(x)^{-2}u\overline v
    \Big)\,\mathrm dx\\
    &=\int_{\Omega}f\overline v\,\mathrm dx,
    \qquad v\in H^1_0(\Omega),
    \end{aligned}
\end{equation}
where the coefficients $D$ and $\beta$ are defined by
\begin{equation*}
  D:=\operatorname{diag}\!\left(
  \gamma_1(x_1)^{-2},\ldots,\gamma_d(x_d)^{-2}\right),
  \qquad
  \beta:=\left(
  \frac{\gamma_1'(x_1)}{\gamma_1(x_1)^3},\ldots,
  \frac{\gamma_d'(x_d)}{\gamma_d(x_d)^3}
  \right)^{\!\top}.
\end{equation*}
The smoothness of the PML profiles implies that the coefficients $D$ and
$\beta$ are bounded. 

For a measurable set $G\subset\mathbb R^d$, we denote the $L^2(G)$ norm by
$\|\cdot\|_{L^2(G)}$ and the complex $L^2(G)$ inner product by
$(\cdot,\cdot)$.
For nonnegative quantities $X$ and $Y$, the notation $X\lesssim Y$ means that
$X\leq CY$, where $C>0$ is independent of $k$ but may depend on the fixed {domain}
geometry, coefficients, and PML data. The value of $C$ may differ between
occurrences. 

\subsection{The RAS--PML method}
The RAS--PML method was introduced and analyzed in
\cite{galkowski2024convergence,galkowski2024schwarz}.
We briefly summarize the main components of RAS--PML below.

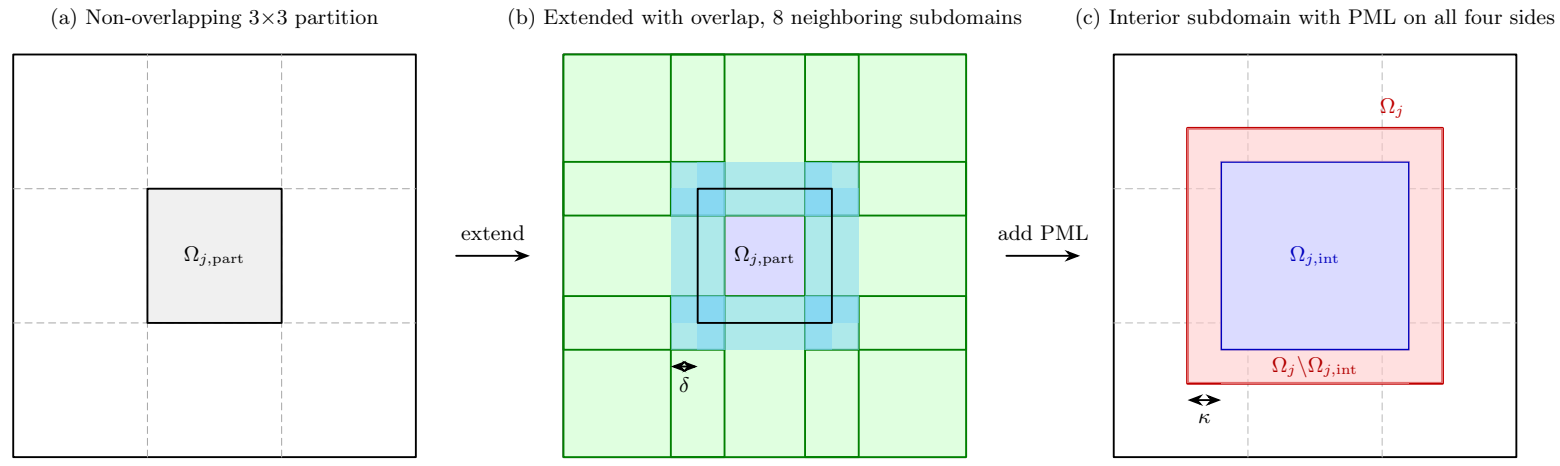
\begin{figure}[htbp!]
    \centering
    \resizebox{0.8\textwidth}{!}{
        \begin{tikzpicture}[font=\small,>=Stealth,line cap=round,line join=round]

\begin{scope}[shift={(0,0)}]
  \draw[thick] (0,0) rectangle (6,6);
  \foreach \x in {2,4} {\draw[densely dashed,gray!70] (\x,0)--(\x,6);}
  \foreach \y in {2,4} {\draw[densely dashed,gray!70] (0,\y)--(6,\y);}
  \fill[gray!12] (2,2) rectangle (4,4);
  \draw[thick] (2,2) rectangle (4,4);
  \node at (3,3) {$\Omega_{j,\mathrm{part}}$};
  \node[above] at (3,6.25) {(a) Non-overlapping 3$\times$3 partition};
\end{scope}

\draw[->,thick] (6.6,3) -- (7.7,3);
\node at (7.15,3.35) {extend};

\begin{scope}[shift={(8.2,0)}]
  \draw[thick] (0,0) rectangle (6,6);
  \foreach \x in {2,4} {\draw[densely dashed,gray!55] (\x,0)--(\x,6);}
  \foreach \y in {2,4} {\draw[densely dashed,gray!55] (0,\y)--(6,\y);}
  \fill[blue!14] (1.6,1.6) rectangle (4.4,4.4);
  \draw[blue!75!black,thick] (1.6,1.6) rectangle (4.4,4.4);
  \node[blue!75!black] at (3,4.75) {$\Omega_{j,\mathrm{int}}$};

  \fill[green!12] (0.0,1.6) rectangle (2.4,4.4);
  \fill[green!12] (3.6,1.6) rectangle (6.0,4.4);
  \fill[green!12] (1.6,0.0) rectangle (4.4,2.4);
  \fill[green!12] (1.6,3.6) rectangle (4.4,6.0);
  \fill[green!12] (0.0,0.0) rectangle (2.4,2.4);
  \fill[green!12] (3.6,0.0) rectangle (6.0,2.4);
  \fill[green!12] (0.0,3.6) rectangle (2.4,6.0);
  \fill[green!12] (3.6,3.6) rectangle (6.0,6.0);
  \draw[green!50!black,thick] (0.0,1.6) rectangle (2.4,4.4);
  \draw[green!50!black,thick] (3.6,1.6) rectangle (6.0,4.4);
  \draw[green!50!black,thick] (1.6,0.0) rectangle (4.4,2.4);
  \draw[green!50!black,thick] (1.6,3.6) rectangle (4.4,6.0);
  \draw[green!50!black,thick] (0.0,0.0) rectangle (2.4,2.4);
  \draw[green!50!black,thick] (3.6,0.0) rectangle (6.0,2.4);
  \draw[green!50!black,thick] (0.0,3.6) rectangle (2.4,6.0);
  \draw[green!50!black,thick] (3.6,3.6) rectangle (6.0,6.0);

  \fill[cyan!45,opacity=0.5] (1.6,2.0) rectangle (2.4,4.0);
  \fill[cyan!45,opacity=0.5] (3.6,2.0) rectangle (4.4,4.0);
  \fill[cyan!45,opacity=0.5] (2.0,1.6) rectangle (4.0,2.4);
  \fill[cyan!45,opacity=0.5] (2.0,3.6) rectangle (4.0,4.4);
  \fill[cyan!45,opacity=0.5] (1.6,1.6) rectangle (2.4,2.4);
  \fill[cyan!45,opacity=0.5] (3.6,1.6) rectangle (4.4,2.4);
  \fill[cyan!45,opacity=0.5] (1.6,3.6) rectangle (2.4,4.4);
  \fill[cyan!45,opacity=0.5] (3.6,3.6) rectangle (4.4,4.4);
  \draw[black,thick] (2.0,2.0) rectangle (4.0,4.0);
  \node at (3,3) {$\Omega_{j,\mathrm{part}}$};
  \draw[<->,thick] (1.6,1.35) -- (2.0,1.35);
  \node at (1.8,1.08) {$\delta$};
  \node[above] at (3,6.25) {(b) Extended with overlap, 8 neighboring subdomains};
\end{scope}

\draw[->,thick] (14.8,3) -- (15.9,3);
\node at (15.35,3.35) {add PML};

\begin{scope}[shift={(16.4,0)}]
  \draw[thick] (0,0) rectangle (6,6);
  \foreach \x in {2,4} {\draw[densely dashed,gray!55] (\x,0)--(\x,6);}
  \foreach \y in {2,4} {\draw[densely dashed,gray!55] (0,\y)--(6,\y);}

  \fill[blue!14] (1.6,1.6) rectangle (4.4,4.4);
  \draw[blue!75!black,thick] (1.6,1.6) rectangle (4.4,4.4);
  \node[blue!75!black] at (3,3) {$\Omega_{j,\mathrm{int}}$};

  \draw[red!75!black,very thick] (1.1,1.1) rectangle (4.9,4.9);
  \node[red!75!black] at (4.15,5.2) {$\Omega_j$};

  \fill[red!22,opacity=0.55] (1.1,1.1) rectangle (1.6,4.9);
  \fill[red!22,opacity=0.55] (4.4,1.1) rectangle (4.9,4.9);
  \fill[red!22,opacity=0.55] (1.6,1.1) rectangle (4.4,1.6);
  \fill[red!22,opacity=0.55] (1.6,4.4) rectangle (4.4,4.9);
  \node[red!70!black] at (3,1.35) {$\Omega_j{\setminus}\Omega_{j,\mathrm{int}}$};
  \draw[<->,thick] (1.1,0.85) -- (1.6,0.85);
  \node at (1.35,0.58) {$\kappa$};

  \node[above] at (3,6.25) {(c) Interior subdomain with PML on all four sides};
\end{scope}

\end{tikzpicture}
    }
    \caption{Construction of the Cartesian covering:
    a non-overlapping partition $\Omega_{j,\mathrm{part}}$ is first
    enlarged to an overlapping subdomain $\Omega_{j,\mathrm{int}}$
    with overlap extension $\delta$, and artificial subdomain boundaries
    are then extended by a local PML of thickness $\kappa$ to form
    $\Omega_j$.}
    \label{fig:cartesian_cover}
\end{figure}

\paragraph{Cartesian covering}
We first decompose the computational domain $\Omega$ into a
non-overlapping Cartesian partition
\[
\overline{\Omega}
=
\bigcup_{j=1}^N \overline{\Omega}_{j,\rm part},
\qquad
    \Omega_{j,\rm part}:=\prod_{i=1}^d
    (a_{i,j,\rm part},b_{i,j,\rm part}),~~j=1,\ldots,N,
\]
where the interiors of distinct subdomains are disjoint.
Each $\Omega_{j,\rm part}$ is then enlarged to an overlapping Cartesian
subdomain
$
\Omega_{j,\rm int}
:=
\prod_{i=1}^d (a_{i,j},b_{i,j})
$ 
by an overlap extension $\delta$ on each artificial side, so that each
subdomain overlaps only with its nearest neighbours; see
Figure~\ref{fig:cartesian_cover}.
Finally, each artificial boundary
$\partial\Omega_{j,\rm int}\not\subset\partial\Omega$
is extended by a Cartesian PML, giving
\[
\Omega_j
=
\prod_{i=1}^d
\big(
a_{i,j}-\kappa_{i,j,l},\,
b_{i,j}+\kappa_{i,j,u}
\big), ~j=1,\dots,N,
\]
where
\[
\kappa_{i,j,l}
=
\begin{cases}
0, & a_{i,j}=a_i-\kappa_{\rm g},\\
\kappa, & a_{i,j}>a_i-\kappa_{\rm g},
\end{cases}
\qquad
\kappa_{i,j,u}
=
\begin{cases}
0, & b_{i,j}=b_i+\kappa_{\rm g},\\
\kappa, & b_{i,j}<b_i+\kappa_{\rm g}.
\end{cases}
\]
Thus, no additional local PML is introduced on portions of
$\partial\Omega_{j,\rm int}$ that already coincide with the outer
boundary of the global PML domain.

On each $\Omega_j$, the local PML profile is defined by
\begin{equation*}
    g_{i,j}(x_i) = \begin{cases}
        g(x_i-b_{i,j}), & x_i \geq b_{i,j}, \\
        0, & x_i \in (a_{i,j}, b_{i,j}), \\
        -g(a_{i,j}-x_i), & x_i \leq a_{i,j},
    \end{cases}
\quad i=1,\dots,d, \quad j=1,\dots,N.
\end{equation*}
On the portion of a boundary subdomain inherited from the global PML, we instead set $g_{i,j}=g_i$.
The local sesquilinear form $a_j(\cdot,\cdot)$ is defined as
in~\eqref{eq:weak}, with integration over $\Omega_j$ and with $D$ and
$\beta$ replaced by the local coefficients $D_j$ and $\beta_j$.
{Since each local profile is either a translation of $g$ or an
inherited global profile, the local coefficients $D_j$ and $\beta_j$
satisfy bounds uniform in $j$.}
\begingroup
For the continuity of $a(\cdot,\cdot)$ in \ref{lem:stable-splitting}, assume
\begin{equation}\label{eq:pml-uniform-coefficient-bound}
  \|D_j\|_{L^\infty(\Omega_j)}
  +k^{-1}\|\beta_j\|_{L^\infty(\Omega_j)}
  +\|c^{-2}\|_{L^\infty(\Omega_j)}\lesssim 1
\end{equation}
holds. The profile used in the numerical experiments of Section~\ref{sec:experiments} satisfies this assumption.
\endgroup

To combine local solutions into a global approximation, we introduce a set of non-negative partition of unity (PoU) functions $\{\chi_j\}$ on $\Omega$ based on the cover $\{\Omega_j\}$, such that $\sum_{j=1}^N \chi_j \equiv 1$, and ${\rm supp}(\chi_j) \subset \Omega_{j, \rm int}$.
Therefore, {$\chi_j \equiv 1$ on the non-overlapping part of $\Omega_j$, namely} $\Omega_{j,\rm novlp} := \Omega_{j,\rm int} \backslash (\bigcup_{l\neq j} \Omega_{l,\rm int})$.

\noindent
\paragraph{Parallel overlapping domain decomposition}

{Starting from $u^0\in H^1_0(\Omega)$, each iteration $n\geq0$
computes $\mathfrak c_j^{n+1}\in H^1_0(\Omega_j)$ in parallel for
$j=1,\ldots,N$:
\begin{equation}\label{eq:helm-pml-local-ras}
  a_j(\mathfrak c_j^{n+1},v_j)
  =(f,v_j)-a(u^n,v_j),
  \qquad v_j\in H^1_0(\Omega_j),
\end{equation}
where the local test functions are extended by zero when used in the global form $a$. The global iterate is then updated by
\begin{equation}\label{eq:update-ras}
  u^{n+1}=u^n+\sum_{j=1}^N\chi_j\mathfrak c_j^{n+1}.
\end{equation}
}

{The following convergence result for the
iteration~\eqref{eq:helm-pml-local-ras}--\eqref{eq:update-ras}
is proved in~\cite{galkowski2024convergence}.}

\begin{theorem}[{Convergence of parallel overlapping methods with PML}]\label{thm:galkowski}
  Assume that the Helmholtz problem is nontrapping and 
  that the fixed
  Cartesian decomposition, PML profiles and smooth PoU
  satisfy the hypotheses of~\cite[Section~1.4 and
  Theorem~1.7]{galkowski2024convergence}. 
  In particular, the overlap and PML
  widths are independent of $k$,  and the local problems are solved exactly.
There is an integer $\mathcal N\geq1$, depending on the decomposition
geometry and ray propagation but independent of $k$ and $f$, such that,
for every $M>0$ and integer $s\geq1$, there exist $C_{M,s}>0$ and $k_0>0$
with
\begin{equation}\label{eq:conv}
  \|u-u^{\mathcal N}\|_{H^s_k(\Omega)}
  \leq C_{M,s}k^{-M}\|u-u^0\|_{H^1_k(\Omega)},
  \qquad k\geq k_0.
\end{equation}
Here $ \|v\|_{H^s_k(\Omega)}^2
  :=\sum_{|\alpha|\leq s}
  \|(k^{-1}\partial)^\alpha v\|_{L^2(\Omega)}^2,~ s\in\mathbb N,~ s\geq1.$
\end{theorem}

In~\cite{galkowski2024convergence}, the additional condition
$g''(x)=0$ for $x\geq\kappa_{\rm lin}$, with
$\kappa_{\rm lin}<\kappa_{\rm g}$, is imposed to make the PML scaling linear
near the outer boundary and thereby enable a reflection argument in the
propagation-of-singularities proof. As explained in
\cite[Remark~2.5]{galkowski2024convergence}, this is a technical restriction
associated with the currently available boundary propagation results, we
therefore omit it in our numerical experiments.

Estimate~\eqref{eq:conv} gives superalgebraic decay in $k$ after a
geometry-dependent number of iterations. For example, when $c\equiv1$ in
two dimensions and the fixed decomposition is an $N_1\times N_2$
checkerboard satisfying~\cite[Theorem~1.4]{galkowski2024convergence},
one may take $\mathcal N=N_1+N_2-1.$ 

\paragraph{{RAS-PML iteration}}
{We briefly recall the discrete RAS--PML formulation; see
\cite[Section~2]{xie_dd29_raspml} for further details.
}
{Let $V_h\subset H^1_0(\Omega)$ be the conforming finite element space
on a uniform mesh $\mathcal T_h$ resolving all subdomains. Set
$V_{h,j}:=\{v_h|_{\Omega_j}:v_h\in V_h\}\cap H^1_0(\Omega_j)$.}
{Let $\mathcal R_{h,j}^{\top}:V_{h,j}\to V_h$ denote nodal extension
by zero, and let $\mathcal R_{h,j}:V_h'\to V_{h,j}'$ be its dual restriction.}
{The weighted extension is defined by}
\begin{equation}
    \widetilde{\mathcal R}_{h,j}^{\top}v_{h,j}
    :=
    \mathcal R_{h,j}^{\top}(\Pi_h(\chi_j v_{h,j})),
    \qquad
    v_{h,j}\in V_{h,j},
    \label{eq:weighted-prolongation}
\end{equation}
{Here $\Pi_h$ is nodal interpolation onto $V_{h,j}$ on the inherited mesh.}
{Let $A_h:V_h\to V_h'$ and $A_{h,j}:V_{h,j}\to V_{h,j}'$ be the
operators induced by $a$ and $a_j$, respectively.
Define $F_h\in V_h'$ by $\langle F_h,v_h\rangle=(f,v_h)$
for all $v_h\in V_h$.}
{Starting from $u_h^0\in V_h$, the RAS--PML iterates satisfy,
for $n\geq0$,}
\begin{equation*}
\begin{aligned}
  &{u_h^{n+1}=u_h^n+\B_h\bigl(F_h-A_hu_h^n\bigr),}
\end{aligned}
\end{equation*}
{with the preconditioner}
\begin{equation}
    \B_h = \sum_{j=1}^N \widetilde{\mathcal R}_{h,j}^{\top} A_{h,j}^{-1} \mathcal R_{h,j}.
    \label{eq:ras-preconditioner}
\end{equation}

{For comparison in the numerical experiments, we also consider RAS-Imp, the corresponding RAS method with impedance conditions on artificial subdomain boundaries in place of local PMLs, see \cite{gong2022convergence,gong2021domain} for further details.}

\section{\texorpdfstring{{Nested RAS--PML method and its analysis}}{Nested RAS--PML method and its analysis}}
\label{sec:nested}

\subsection{\texorpdfstring{Nested RAS--PML method}{Nested RAS--PML method}}
\label{ssec:nested-method}
{Theorem~\ref{thm:galkowski} provides the starting point for our
construction. To reduce the cost of large local solves, we approximate
them recursively using RAS--PML.}

The global PML problem forms the root of the hierarchy at level $0$. For
clarity, we describe a two-dimensional $2\times2$ construction; the method
extends directly to other dimensions and nesting patterns. Each
non-terminal parent at level $\ell$, $0\leq\ell<L$, is covered by four
overlapping Cartesian child subdomains. RAS--PML iterations on this child
covering approximate the parent solve, with all child problems solved
concurrently at each iteration. The recursion terminates at level $L$,
where the local problems are solved exactly.

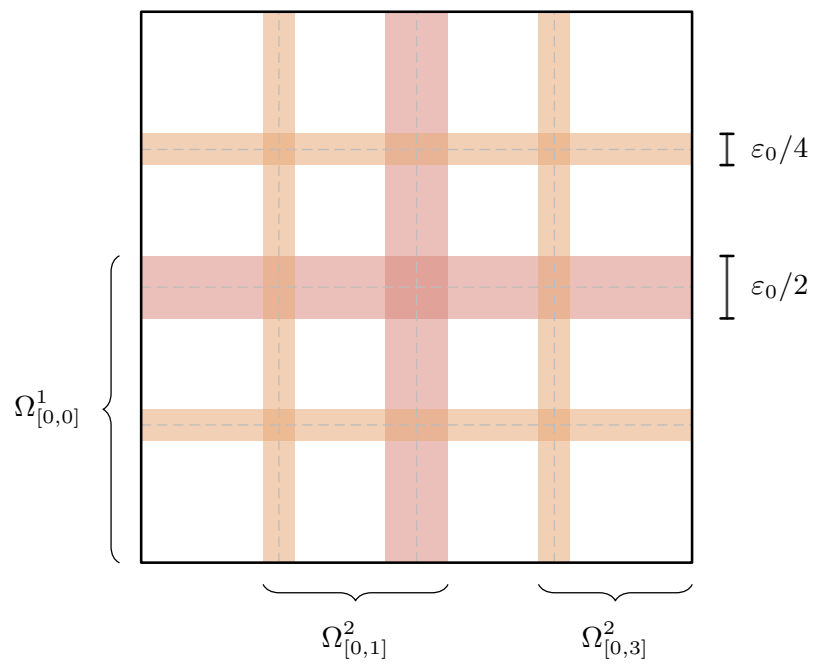
\begin{figure}[H]
    \centering
    \resizebox{0.42\textwidth}{!}{%
\begin{tikzpicture}[font=\small,line cap=round,line join=round]

\def\W{6}
\pgfmathsetmacro{\hCoarse}{0.34}
\pgfmathsetmacro{\hFine}{0.17}
\def\ticklen{0.07}

\definecolor{coarseband}{RGB}{218,138,128}
\definecolor{fineband}{RGB}{232,168,118}

\fill[coarseband,opacity=0.54] ({3-\hCoarse},0) rectangle ({3+\hCoarse},\W);
\fill[coarseband,opacity=0.54] (0,{3-\hCoarse}) rectangle (\W,{3+\hCoarse});
\foreach \xc in {1.5,4.5} {%
  \fill[fineband,opacity=0.54] ({\xc-\hFine},0) rectangle ({\xc+\hFine},\W);
}
\foreach \yc in {1.5,4.5} {%
  \fill[fineband,opacity=0.54] (0,{\yc-\hFine}) rectangle (\W,{\yc+\hFine});
}

\foreach \x in {1.5,3,4.5} {\draw[densely dashed,gray!52] (\x,0) -- (\x,\W);}
\foreach \y in {1.5,3,4.5} {\draw[densely dashed,gray!52] (0,\y) -- (\W,\y);}
\draw[thick] (0,0) rectangle (\W,\W);

\draw[decorate,decoration={brace,amplitude=5pt}]
  (-0.24,0) -- (-0.24,{3+\hCoarse})
  node[midway,anchor=east,xshift=-8pt] {$\Omega_{[0,0]}^1$};
\draw[decorate,decoration={brace,mirror,amplitude=5pt}]
  ({1.5-\hFine},-0.24) -- ({3+\hCoarse},-0.24)
  node[midway,anchor=north,yshift=-8pt] {$\Omega_{[0,1]}^2$};
\draw[decorate,decoration={brace,mirror,amplitude=5pt}]
  ({4.5-\hFine},-0.24) -- (6,-0.24)
  node[midway,anchor=north,yshift=-8pt] {$\Omega_{[0,3]}^2$};

\coordinate (evA) at ({\W+0.38},{3-\hCoarse});
\coordinate (evB) at ({\W+0.38},{3+\hCoarse});
\draw[thick,black!75] (evA) -- (evB);
\draw[thick] ($(evA)+(-\ticklen,0)$) -- ($(evA)+(\ticklen,0)$);
\draw[thick] ($(evB)+(-\ticklen,0)$) -- ($(evB)+(\ticklen,0)$);
\node[right=4pt] at ($(evA)!0.5!(evB)$) {$\varepsilon_0/2$};
\coordinate (ehA) at ({\W+0.38},{4.5-\hFine});
\coordinate (ehB) at ({\W+0.38},{4.5+\hFine});
\draw[thick,black!75] (ehA) -- (ehB);
\draw[thick] ($(ehA)+(-\ticklen,0)$) -- ($(ehA)+(\ticklen,0)$);
\draw[thick] ($(ehB)+(-\ticklen,0)$) -- ($(ehB)+(\ticklen,0)$);
\node[right=4pt] at ($(ehA)!0.5!(ehB)$) {$\varepsilon_0/4$};

\end{tikzpicture}%
    }
    \caption{Two levels of the nested Cartesian hierarchy. At each level, the
    widths of the extended regions, consisting of the PML and overlap layers,
    are halved. Here $\varepsilon_0$ denotes the reference width of an
    extended region at level $0$, so its width at level $\ell$ is
    $2^{-\ell}\varepsilon_0$. In particular,
    $\Omega_{[0,1]}^2\subset\Omega_{[0,0]}^1\subset\Omega$.}
    \label{fig:nested_ras}
\end{figure}

{We first introduce the notation needed to define the nested
RAS--PML iteration.} For
$0\leq\ell\leq L$, let
$
  \mathcal I_\ell
  :=\{0,\ldots,2^\ell-1\}\times\{0,\ldots,2^\ell-1\}.
$ 
A vector $\vec J\in\mathcal I_\ell$ is the global Cartesian index of a
subdomain in the level-$\ell$ covering, and the corresponding subdomain is
denoted by $\Omega_{\vec J}^{\ell}$. In particular,
$\mathcal I_0=\{\vec 0\}$ and $\Omega_{\vec 0}^{0}:=\Omega$.
Let $\mathcal S:=\{0,1\}\times\{0,1\}$ denote the local child-index
set. Given a
level-$\ell$ parent $\Omega_{\vec J}^{\ell}$, where $0\leq\ell<L$, each
$\vec j\in\mathcal S$ identifies one of its four children. The global index of
this child in the level-$(\ell+1)$ covering is
$\vec J_{\vec j}:=2\vec J+\vec j,$ 
and the corresponding child subdomain is
$\Omega_{\vec J_{\vec j}}^{\ell+1}$. Each child is equipped with local PML
regions constructed as in the one-level setting above. Thus, $\vec J$ locates
the parent globally within level $\ell$, whereas $\vec j$ selects a child
locally within that parent.
This indexing is illustrated in Figure~\ref{fig:nested_ras}.

For each $\Omega_{\vec J}^{\ell}$, $1\leq\ell\leq L$, let
$a_{\vec J}^{\ell}(\cdot,\cdot)$ denote its local PML sesquilinear form,
and set $A_{h,\vec J}^{\ell}
  :V_{h,\vec J}^{\ell}\rightarrow(V_{h,\vec J}^{\ell})'$
with $
  V_{h,\vec J}^{\ell}
  :=\{v_h|_{\Omega_{\vec J}^{\ell}}:v_h\in V_h\}
  \cap H_0^1(\Omega_{\vec J}^{\ell})
$, such that
\[
\begin{aligned}
  \langle A_{h,\vec J}^{\ell}u_h,v_h\rangle
  =a_{\vec J}^{\ell}(u_h,v_h),
  ~ \forall u_h,v_h\in V_{h,\vec J}^{\ell}.
\end{aligned}
\]
At level $0$, all quantities refer to their global counterparts, and parent
subscripts are omitted.
Let
$\{\chi_{\vec J_{\vec j}}^{\ell+1}\}_{\vec j\in\mathcal S}$ be a partition
of unity subordinate to the child covering of $\Omega_{\vec J}^{\ell}$.
As in the one-level setting, the operator
$\mathcal R_{h,\vec J_{\vec j}}^{\ell+1,\top}:
V_{h,\vec J_{\vec j}}^{\ell+1}\to V_{h,\vec J}^{\ell}$ is nodal extension by zero,
$\mathcal R_{h,\vec J_{\vec j}}^{\ell+1}:
(V_{h,\vec J}^{\ell})'\to(V_{h,\vec J_{\vec j}}^{\ell+1})'$
{is its adjoint restriction}, and
$\widetilde{\mathcal R}_{h,\vec J_{\vec j}}^{\ell+1,\top}$ is the weighted
prolongation defined as in~\eqref{eq:weighted-prolongation} with weight
$\chi_{\vec J_{\vec j}}^{\ell+1}$.

The general nested iteration uses a level-dependent sequence
$\{\mathcal N_\ell\}_{\ell=1}^{L-1}$ of inner iteration counts.
We denote by Nested($q$) its fixed-count specialization
$\mathcal N_\ell=q$ for all $1\leq\ell<L$, while the level-$0$ iteration continues until the stopping
criterion is met.
Algorithm~\ref{alg:nestedras} summarizes the recursive
procedure, where $\ell$ denotes the level of the current parent problem.
\SetKwProg{Fn}{Function}{:}{}
\begin{algorithm2e}[h!]
  \small
    \caption{Nested RAS--PML iteration for a $2\times2$ Cartesian hierarchy}
    \label{alg:nestedras}
    \KwIn{Initial iterate $u_h^0$, source term $f$, domain $\Omega$, overlap extension $\delta$, PML parameters, hierarchy depth $L$, levelwise inner iteration counts
$\{\mathcal N_\ell\}_{\ell=1}^{L-1}$,  stopping tolerance \texttt{tol}, maximum number of iterations \texttt{maxit}}
    \KwOut{Approximate solution $u_h$}
    {
    \Fn{\textnormal{\textsc{NestedRAS}}$(u_{\vec J}^{\ell,0},
    g_{\vec J}^{\ell},\ell)$}{
    Set the child indices $\vec J_{\vec j}$ as above for
    $\vec j\in\mathcal S$\;
    \For{$n=0,1,2,\ldots$}{
        \ForEach{$\vec j\in\mathcal S$ \Comment{Parallel child solves}}{
            Set $g_{\vec J_{\vec j}}^{\ell+1,n}
            \gets\mathcal R_{h,\vec J_{\vec j}}^{\ell+1}
            (g_{\vec J}^{\ell}-A_{h,\vec J}^{\ell}
            u_{\vec J}^{\ell,n})$\;
            \eIf{$\ell+1=L$ \Comment{Terminal level}}{
                $\mathfrak c_{\vec J_{\vec j}}^{\ell+1,n+1}
                \gets(A_{h,\vec J_{\vec j}}^{\ell+1})^{-1}
                g_{\vec J_{\vec j}}^{\ell+1,n}$\;
            }{
                $\mathfrak c_{\vec J_{\vec j}}^{\ell+1,n+1}
                \gets$
                \Call{NestedRAS}
                {$0,g_{\vec J_{\vec j}}^{\ell+1,n},\ell+1$}\;
            }
        }
        $u_{\vec J}^{\ell,n+1}\gets u_{\vec J}^{\ell,n}
        +\displaystyle\sum_{\vec j\in\mathcal S}
        \widetilde{\mathcal R}_{h,\vec J_{\vec j}}^{\ell+1,\top}
        \mathfrak c_{\vec J_{\vec j}}^{\ell+1,n+1}$\;
        \If{{($\ell=0$ and (relative residual $<\texttt{tol}$ or
        $n+1=\texttt{maxit}$)) or
        ($\ell\geq1$ and $n+1=\mathcal N_\ell$)}}{
            break\;
        }
    }
    {\Return $u_{\vec J}^{\ell,n+1}$\;}
    }
    {\Return \Call{NestedRAS}{$u_h^0,F_h,0$}\;}
    }
\end{algorithm2e}

The operator notation below is the algebraic form of this recursion. At the
terminal level, $(\widetilde A_{h,\vec J}^{L})^{-1}
  :=(A_{h,\vec J}^{L})^{-1}.$ 
For a level-$\ell$ parent, $0\leq\ell<L$, define the exact and nested
RAS--PML preconditioners by
\begin{equation}\label{eq:levelwise-ras-preconditioners}
\begin{aligned}
  \B_{h,\vec J}^{\ell}
  &:=\sum_{\vec j\in\mathcal S}
  \widetilde{\mathcal R}_{h,\vec J_{\vec j}}^{\ell+1,\top}
  (A_{h,\vec J_{\vec j}}^{\ell+1})^{-1}
  \mathcal R_{h,\vec J_{\vec j}}^{\ell+1},\\
  \Bt_{h,\vec J}^{\ell}
  &:=\sum_{\vec j\in\mathcal S}
  \widetilde{\mathcal R}_{h,\vec J_{\vec j}}^{\ell+1,\top}
  (\widetilde A_{h,\vec J_{\vec j}}^{\ell+1})^{-1}
  \mathcal R_{h,\vec J_{\vec j}}^{\ell+1}.
\end{aligned}
\end{equation}
Here $A_{h,\vec J}^{\ell}$ is the level-$\ell$ discrete operator. The exact
RAS--PML preconditioner $\B_{h,\vec J}^{\ell}$ is assembled from the exact
child solvers $(A_{h,\vec J_{\vec j}}^{\ell+1})^{-1}$, whereas the nested
RAS--PML preconditioner $\Bt_{h,\vec J}^{\ell}$ is assembled from the
recursively constructed child solvers
$(\widetilde A_{h,\vec J_{\vec j}}^{\ell+1})^{-1}$. For
$\ell=L-1,\ldots,1$, define the recursive level-$\ell$ solver (equivalently,
the approximate inverse) $(\widetilde A_{h,\vec J}^{\ell})^{-1}$ through its
action: $(\widetilde A_{h,\vec J}^{\ell})^{-1}g_{\vec J}^{\ell}$ is the output
of $\mathcal N_\ell$ Richardson iterations with preconditioner
$\Bt_{h,\vec J}^{\ell}$ and zero initial guess. Its error-propagation operator
therefore satisfies
\begin{equation}\label{eq:recursive-approximate-inverse}
  I-(\widetilde A_{h,\vec J}^{\ell})^{-1}A_{h,\vec J}^{\ell}
  =\left(I-\Bt_{h,\vec J}^{\ell}A_{h,\vec J}^{\ell}\right)^{\mathcal N_\ell}.
\end{equation}
At the root level $\ell=0$, we write $ \B_h:=\B_{h,\vec 0}^{0},~\Bt_h:=\Bt_{h,\vec 0}^{0}.$

\subsection{\texorpdfstring{Main results}{Main results}} 
\label{ssec:main-results}
{We now state the two main results for the nested RAS--PML method. The first result proves that Nested(1) is equivalent to applying the RAS-PML algorithm to the finest level domain decomposition, with a certain (non-standard) PoU.} The second controls convergence under sufficiently accurate recursive local solution.

\subsubsection{\texorpdfstring{{Algebraic equivalence}}{Algebraic equivalence}}

Let $\boldsymbol j=(\vec j_1,\ldots,\vec j_L)\in\mathcal S^L$ denote a
branch of the hierarchy, with
$\vec J_1=\vec j_1$ and
$\vec J_\ell=2\vec J_{\ell-1}+\vec j_\ell$ for $2\leq\ell\leq L$.
Define the composite restriction and weighted extension by
\[
  \widehat{\mathcal R}_{h,\boldsymbol j}
  :=\mathcal R_{h,\vec J_L}^{L}\cdots
  \mathcal R_{h,\vec J_1}^{1},
  \qquad
  \widehat{\widetilde{\mathcal R}}_{h,\boldsymbol j}^{\top}
  :=\widetilde{\mathcal R}_{h,\vec J_1}^{1,\top}\cdots
  \widetilde{\mathcal R}_{h,\vec J_L}^{L,\top}.
\]
{The following theorem shows that Nested(1) can be implemented
as a one-level RAS--PML method without recursive transfers, using
the local problems and composite transfer operators induced by
the hierarchy. The proof appears in Subsection~\ref{ssec:equivalence-proof}.}
\begin{theorem}[Equivalence of Nested(1) and one-level RAS-PML]
\label{thm:nested-ras-equivalence}
Suppose that one-level RAS-PML uses the finest-level local problems of
Nested(1), together with the composite restriction
$\widehat{\mathcal R}_{h,\boldsymbol j}$ and weighted extension
$\widehat{\widetilde{\mathcal R}}_{h,\boldsymbol j}^{\top}$ induced by the
hierarchy. If the two methods start from the same initial iterate, then, for
every $n\geq 0$, the approximate solution
$u_{h,\mathrm{Nested}(1)}^n$ produced by Nested(1) coincides with the
approximate solution $u_{h,\mathrm{RAS}}^n$ produced by the one-level RAS-PML
method.
\end{theorem}

\subsubsection{\texorpdfstring{{Convergence with {sufficiently} accurate recursive local solves}}{Convergence with accurate recursive local solves}}

{Under suitable stability and contraction assumptions,
Theorem~\ref{thm:nested-conv} establishes convergence of the nested iteration when the recursive local solves are sufficiently accurate.}

{In the normalized dyadic hierarchy, the nominal side length at
level $\ell$ is $2^{-\ell}$. Rescaling this length to one gives the
effective wavenumber $k2^{-\ell}$. The high-frequency stability and
contraction estimates require this parameter to be sufficiently large.
We choose $K_*>0$, independent of $k$, sufficiently large for
the superalgebraic estimate in Theorem~\ref{thm:galkowski}
to yield contraction. Uniform discrete contraction across
levels is assumed below.} For $k\geq2K_*$, choose the hierarchy depth as
\begin{equation}\label{eq:hierarchy-depth}
  L=L(k):=\left\lfloor\log_2\!\left(\frac{k}{K_*}\right)\right\rfloor,
  \qquad
  K_*\leq\frac{k}{2^L}<2K_*.
\end{equation}
Then $L=\mathcal O(\log k)$, and every non-terminal level satisfies
$  \dfrac{k}{2^\ell}\geq2K_*,~ 0\leq\ell<L.$

\begin{assumption}[Uniform discrete high-frequency contraction]
\label{ass:uniform-discrete-contraction}
There exist $\mathcal N\in\mathbb N$ and $\rho\in(0,1)$, independent of
$k$, $L$, and the level, such that
\[
  \max_{0\leq\ell<L}
  \max_{\vec J\in\mathcal I_\ell}
  \left\|
    \left(I-\B_{h,\vec J}^{\ell}
    A_{h,\vec J}^{\ell}\right)^{\mathcal N}
  \right\|_{
    H^1_k(\Omega_{\vec J}^{\ell})\to
    H^1_k(\Omega_{\vec J}^{\ell})}
  \leq\rho.
\]
\end{assumption}

\begingroup
\begin{remark}[Scope of Assumption~\ref{ass:uniform-discrete-contraction}]
\label{rmk:high-frequency-contraction}
For a fixed decomposition, Theorem~\ref{thm:galkowski} gives contraction
of the continuous iteration after a fixed number of steps for sufficiently
large $k$.
Here the hierarchy depth grows with $k$ as specified in~\eqref{eq:hierarchy-depth}.
These estimates do not by themselves provide a uniform contraction bound
for this growing discrete hierarchy.
We therefore assume uniform power contractivity across all non-terminal
levels, as stated in Assumption~\ref{ass:uniform-discrete-contraction}.
\end{remark}
\endgroup

\begin{theorem}[Convergence of the nested RAS--PML iteration]
\label{thm:nested-conv}
Suppose that Assumption~\ref{ass:uniform-discrete-contraction} holds, {\eqref{eq:pml-uniform-coefficient-bound} holds uniformly for the local problems on each parent subdomain and its child covering,
and that the discrete inf--sup constant $\mu_{h,\vec J}^{\ell}$ of the local problem on
$\Omega_{\vec J}^{\ell}$ satisfies}
\[
  {(\mu_{h,\vec J}^{\ell})^{-1}
  \lesssim \frac{k}{2^\ell},}
  \qquad
  \vec J\in\mathcal I_\ell,
  \quad 1\leq\ell<L,
\]
with a uniform hidden constant. Choose $q_1>0$, independently of $k$ and $L$,
such that $q_1\log\left(\frac{2}{1+\rho}\right)>\mathcal N.$ 
Then, for all sufficiently large $q_0>0$ independent of $k$ and $L$, define
the iteration counts by
\begin{equation}\label{eq:inner-iteration-count}
  \mathcal N_\ell
  :=\mathcal N\left\lceil
  q_0+q_1\log\left(1+\frac{k}{2^\ell}\right)
  \right\rceil,
  \qquad 1\leq\ell<L.
\end{equation}
The global nested iteration then satisfies
\[
  \left\|\left(I-\Bt_hA_h\right)^{\mathcal N}\right\|
  \leq\frac{1+\rho}{2}<1.
\]
\end{theorem}
Theorem~\ref{thm:nested-conv} gives uniform contraction in blocks of
$\mathcal N$ outer iterations for sufficiently accurate recursive local
solves. {The hierarchy depth $L$ in Algorithm~\ref{alg:nestedras}}
grows logarithmically with $k$, and the sufficient
inner iteration counts depend on the effective wavenumber at each level.
Its proof is given in Subsection~\ref{ssec:conv}.

\subsection{\texorpdfstring{{Proof of algebraic equivalence}}{Proof of algebraic equivalence}}
\label{ssec:equivalence-proof}

\begingroup
\begin{proof}[Proof of Theorem~\ref{thm:nested-ras-equivalence}]
Along a branch, the parent and child spaces use the inherited global
mesh. Nodal interpolation therefore gives, for
$v\in V_{h,\vec J_L}^{L}$,
\[
  \widehat{\widetilde{\mathcal R}}_{h,\boldsymbol j}^{\top}v
  =
  \sum_{x\in\mathcal V_{h,\vec J_L}^{L}}
  \left(\prod_{\ell=1}^L\chi_{\vec J_\ell}^{\ell}(x)\right)
  v(x)\varphi_x.
\]
Here $\mathcal V_{h,\vec J_L}^{L}$ is the terminal nodal set and
$\varphi_x$ is the global nodal basis function. This follows by
applying~\eqref{eq:weighted-prolongation} successively: at each retained
node, a weighted extension multiplies the current nodal value by the
parent-level PoU weight. Extending child weights by zero outside their
parent and summing successively from the deepest level gives
  $\sum_{\boldsymbol j\in\mathcal S^L}
  \prod_{\ell=1}^L\chi_{\vec J_\ell}^{\ell}=1.$ 
The composite weights thus form a nonnegative PoU subordinate to the
finest covering.

Suppose that both methods have iterate $u_h^n$ and set
$r_h^n:=F_h-A_hu_h^n$. Every non-terminal inner call of Nested(1)
starts from zero and performs one iteration. Its initial residual is
its right-hand side, so the right-hand side reaching a terminal problem
is $\widehat{\mathcal R}_{h,\boldsymbol j}r_h^n$.
Expanding the recursive corrections, with exact terminal solves, yields
\begin{equation*}
  u_{h,\mathrm{Nested}(1)}^{n+1}
  =u_h^n+
  \sum_{\boldsymbol j\in\mathcal S^L}
  \widehat{\widetilde{\mathcal R}}_{h,\boldsymbol j}^{\top}
  (A_{h,\vec J_L}^{L})^{-1}
  \widehat{\mathcal R}_{h,\boldsymbol j}
  (F_h-A_hu_h^n).
\end{equation*}
This is exactly the one-level update with the stated local problems
and composite transfer operators. {Both methods produce the same next iterate from the same current
iterate, including when $L=1$.
Equality of the initial iterates therefore proves the result.}
\end{proof}
\endgroup


\subsection{Proof of convergence of the nested RAS--PML iteration}
{We prove Theorem~\ref{thm:nested-conv} in three stages.
Section~\ref{ssec:infsup} recalls the continuous inf--sup estimate and
discusses the discrete local stability assumptions under rescaling.
Section~\ref{ssec:bdd} bounds the exact one-level preconditioned operator
using weighted-assembly and local-correction estimates.
Finally, Section~\ref{ssec:conv} estimates the perturbations caused by
inexact child solves and uses downward induction to control the
recursive solver errors, completing the convergence proof.}

\subsubsection{Continuous and discrete inf-sup conditions}
\label{ssec:infsup}
\begin{lemma}[Continuous inf--sup condition]\label{lem:cts_infsup}
Assume that the Helmholtz {problem} is non-trapping {in the sense of \cite{galkowski2024convergence}} and that the fixed
Cartesian PML geometry and smooth scaling profiles satisfy the
assumptions of~\cite[Lemma~4.5]{galkowski2024convergence}. Then, for all
sufficiently large $k$,
\[
  \sup_{0\ne v\in H^1_0(\Omega)}
  \frac{|a(u,v)|}{\|v\|_{H^1_k(\Omega)}}
  \geq\mu\|u\|_{H^1_k(\Omega)},
  \qquad u\in H^1_0(\Omega),
  \qquad \mu\gtrsim k^{-1}.
\]
\end{lemma}

%

{We assume that the global Cartesian PML finite element problem
and all local finite element problems are uniquely solvable.
For completeness, the discrete inf--sup constant of the global problem
is given by
\[
  \mu_h:=\inf_{0\ne u_h\in V_h}\sup_{0\ne v_h\in V_h}
  \frac{|a(u_h,v_h)|}
  {\|u_h\|_{H^1_k(\Omega)}\|v_h\|_{H^1_k(\Omega)}}.
\]
The local constants $\mu_{h,j}$ and $\mu_{h,\vec J}^{\ell}$
are defined analogously.}

%

{Under suitable regularity and approximation conditions, existing analyses
bound the discrete inf--sup constant below by a $k$-independent multiple
of the continuous one; see~\cite[Theorem~4.2]{Melenk_2010}
and~\cite[Section~4.3]{spence2026preconditioning}.
For the Cartesian PML hierarchy considered here, we assume the required
local bounds to hold uniformly.}

\begingroup
\begin{remark}[Local scaling and stability]\label{rmk:subdomain-infsup}
Let $H_j:=\operatorname{diam}(\Omega_j)$, set
$\widehat\Omega_j:=(\Omega_j-x_j)/H_j$, and write
$\widehat u(\widehat x):=u(x_j+H_j\widehat x)$.
The rescaled coefficients are
$\widehat D=D_j(x_j+H_j\widehat x)$,
$\widehat\beta=H_j\beta_j(x_j+H_j\widehat x)$, and
$\widehat c=c(x_j+H_j\widehat x)$.
With effective wavenumber $kH_j$,
\[
  a_j(u,v)=H_j^d\widehat a_j(\widehat u,\widehat v),
  \qquad
  \|u\|_{H^1_k(\Omega_j)}
  =H_j^{d/2}\|\widehat u\|_{H^1_{kH_j}(\widehat\Omega_j)}.
\]

{The continuous inf--sup quotient is unchanged by rescaling.
If the rescaled family satisfies Lemma~\ref{lem:cts_infsup} with uniform
constants, then $\mu_j\gtrsim(kH_j)^{-1}$ for sufficiently large $kH_j$.
For the discrete local problems, the corresponding bound
$\mu_{h,j}^{-1}\lesssim kH_j$ is an additional assumption.
In the normalized dyadic hierarchy with $H_\ell\simeq2^{-\ell}$,  we assume that, for $1\leq\ell<L$,
the following bound holds uniformly in $k$ and the level:
\begin{equation*}
  (\mu_{h,\vec J}^{\ell})^{-1}
  \lesssim kH_\ell\lesssim \frac{k}{2^\ell}.
\end{equation*}
}
\end{remark}
\endgroup

\subsubsection[Boundedness of the exact one-level preconditioned operator]
{Boundedness of the exact one-level {preconditioned operator}}
\label{ssec:bdd}

Recall the exact one-level preconditioner $\B_h$
from~\eqref{eq:ras-preconditioner}. Define the local projection
$Q_{h,j}:V_h\to V_{h,j}$, for $v_h\in V_h$, by
\begin{equation*}
  a_j(Q_{h,j} v_h,w_{h,j})
  =a(v_h,\mathcal R_{h,j}^{\top}w_{h,j}),
  \qquad\forall w_{h,j}\in V_{h,j},
\end{equation*}
and set
$Q_h v_h=\sum_{j=1}^N\widetilde{\mathcal R}_{h,j}^{\top}Q_{h,j} v_h$.
Then
$Q_{h,j} v_h=A_{h,j}^{-1}\mathcal R_{h,j}A_h v_h$ and
$Q_h=\B_hA_h$.
For each element $\tau\in\mathcal T_h$, let $h_\tau:=\operatorname{diam}(\tau)$,
and define
\[
  \mathcal T_{h,j}:=
  \{\tau\in\mathcal T_h:\tau\subset\overline{\Omega_j}\},
  \qquad
  h_j:=\max_{\tau\in\mathcal T_{h,j}}h_\tau,
  \qquad
  h:=\max_{\tau\in\mathcal T_h}h_\tau
     =\max_{1\le j\le N}h_j.
\]
For each $j$, let $\delta_j$
denote the overlap width of $\Omega_j$, and set
$\delta:=\min_{1\le j\le N}\delta_j$.
Let $\{\chi_j\}$ be the smooth PoU
introduced above, satisfying
\[
  \|\nabla\chi_j\|_{L^\infty(\Omega_j)}
  \lesssim\delta_j^{-1},
  \qquad j=1,\ldots,N.
\]
Its action on the discrete spaces is realized through nodal interpolation
as in~\eqref{eq:weighted-prolongation}.

\begin{assumption}\label{ass:hkd}
For every $j$, assume that
\[
  kh_j\lesssim1,
  \qquad k\delta_j\gtrsim1,
\]
with constants independent of $j$. Equivalently, $kh\lesssim1$ and
$k\delta\gtrsim1$.
\end{assumption}

For these PoU functions, we use the following standard stability and
interpolation estimates; see
\cite[Lemmas 4.2 and 4.4]{gong2021domain} and
\cite[Lemma 3.3]{graham2020domain}.
\begin{lemma}[PoU estimates]\label{lem:interpolation}
For all $v_{h,j}\in V_{h,j}$,
\begin{equation*}
  \|\chi_j v_{h,j}\|_{H^1_k(\Omega_j)}
  \lesssim(1+(k\delta_j)^{-1})\|v_{h,j}\|_{H^1_k(\Omega_j)},
\end{equation*}
\[
  \|(I-\Pi_h)(\chi_j v_{h,j})\|_{H^1_k(\Omega_j)}
  \lesssim(1+kh_j)\Bigl(\frac{h_j}{\delta_j}\Bigr)
  \|v_{h,j}\|_{H^1_k(\Omega_j)}.
\]
\end{lemma}

We assume that the covering $\{\Omega_j\}_{j=1}^N$ has uniformly finite
overlap. More precisely, define
\[
  \Lambda(i):=\{j:\Omega_i\cap\Omega_j\neq\emptyset\},
  \qquad
  \Lambda:=\max_{1\le i\le N}\#\Lambda(i)
\]
and suppose that $\Lambda$ is bounded independently of $N$.

\begin{lemma}[Stability of weighted assembly]\label{lem:weighted-assembly}
Under Assumption~\ref{ass:hkd}, every family $w_{h,j}\in V_{h,j}$ satisfies
\[
\begin{aligned}
  \left\|\sum_{j=1}^N
  \widetilde{\mathcal R}_{h,j}^{\top}w_{h,j}\right\|_{H^1_k(\Omega)}^2
  &\lesssim\Lambda\sum_{j=1}^N
  \|w_{h,j}\|_{H^1_k(\Omega_j)}^2.
\end{aligned}
\]
\end{lemma}
{The estimate follows from the definition of the weighted
prolongation, Lemma~\ref{lem:interpolation} and
Assumption~\ref{ass:hkd}, together with the finite-overlap
argument in \cite[Lemma 4.2]{grahamDomainDecompositionPreconditioning2017}.}

\begin{lemma}\label{lem:stable-splitting}

{Assume~\eqref{eq:pml-uniform-coefficient-bound} and its
global counterpart hold with the constant independent of $k$, the hierarchy
depth, level, and subdomain index, and that}
{all local finite element problems are uniquely solvable, with discrete inf-sup constants $\mu_{h,j}>0$.}
Then, for all $v_h\in V_h$,
\[
  \sum_{j=1}^N\|Q_{h,j} v_h\|_{H^1_k(\Omega_j)}^2
  \lesssim\Lambda\max_{1\le j\le N}\mu_{h,j}^{-2}
  \|v_h\|_{H^1_k(\Omega)}^2.
\]
\end{lemma}
{From \eqref{eq:pml-uniform-coefficient-bound}, the Cauchy-Schwarz inequality shows that \(a(\cdot,\cdot)\) is continuous with a constant independent of \(k\). Together with the local discrete inf--sup condition and the finite overlap property, this shows the above result.}
%
{For $Q_h=\B_hA_h$, we can easily get the following result from Lemmas~\ref{lem:weighted-assembly}
and~\ref{lem:stable-splitting}.}
\begin{theorem}[Boundedness of the exact one-level
{preconditioned operator}]\label{thm:bdd}
Under Assumption~\ref{ass:hkd} and the hypotheses of
Lemma~\ref{lem:stable-splitting},
the preconditioned operator \(Q_h=\B_hA_h\) satisfies
\[
  \|\B_hA_h\|\lesssim\Lambda
  \max_{1\le j\le N}\mu_{h,j}^{-1},
  \qquad
  \|I-\B_hA_h\|\lesssim1+\Lambda
  \max_{1\le j\le N}\mu_{h,j}^{-1},
\]
where the operator norms are induced by $\|\cdot\|_{H^1_k(\Omega)}$.
\end{theorem}

\subsubsection{Perturbation estimates and proof of convergence}
\label{ssec:conv}

\begingroup
The proof follows three steps. We first bound the perturbation caused by
inexact child solves, then estimate its effect on a fixed block of
iterations, and finally choose the inner accuracies by downward induction.
\endgroup

For $1\leq\ell\leq L$, denote the maximal local-solver error in
operator norm by
\[
  E_\ell
  :=\max_{\vec J\in\mathcal I_\ell}\left\|
  I-(\widetilde A_{h,\vec J}^{\ell})^{-1}
  A_{h,\vec J}^{\ell}
  \right\|_{H^1_k(\Omega_{\vec J}^{\ell})\to H^1_k(\Omega_{\vec J}^{\ell})}.
\]
We omit the underlying spaces from the
notation for operator norms
whenever they are clear from the context.
\begin{lemma}\label{lem:B-Bh}
Let $0\leq\ell<L-1$ and fix a level-$\ell$ parent
$\Omega_{\vec J}^{\ell}$. Assume that the hypotheses of the levelwise
counterparts of Lemmas~\ref{lem:weighted-assembly}
and~\ref{lem:stable-splitting} hold for its child covering. Then
\[
\begin{aligned}
  &
  \left\|
  (\B_{h,\vec J}^{\ell}-\Bt_{h,\vec J}^{\ell})
  A_{h,\vec J}^{\ell}
  \right\|_{\substack{
    H^1_k(\Omega_{\vec J}^{\ell})\to
    H^1_k(\Omega_{\vec J}^{\ell})}} \\
  &\quad\lesssim\Lambda\max_{\vec j\in\mathcal S}\left[
  (\mu_{h,\vec J_{\vec j}}^{\ell+1})^{-1}
  \left\|
  I-(\widetilde A_{h,\vec J_{\vec j}}^{\ell+1})^{-1}
  A_{h,\vec J_{\vec j}}^{\ell+1}
  \right\|_{\substack{
    H^1_k(\Omega_{\vec J_{\vec j}}^{\ell+1})\to
    H^1_k(\Omega_{\vec J_{\vec j}}^{\ell+1})}}
  \right].
\end{aligned}
\]

\end{lemma}
\begin{proof}
For $v_{h,\vec J}\in V_{h,\vec J}^{\ell}$, the definitions of
$\B_{h,\vec J}^{\ell}$ and $\Bt_{h,\vec J}^{\ell}$
in~\eqref{eq:levelwise-ras-preconditioners} give
\[
\begin{aligned}
  &(\B_{h,\vec J}^{\ell}-\Bt_{h,\vec J}^{\ell})
  A_{h,\vec J}^{\ell}v_{h,\vec J}\\
  &\quad=\sum_{\vec j\in\mathcal S}
  \widetilde{\mathcal R}_{h,\vec J_{\vec j}}^{\ell+1,\top}
  \left[(A_{h,\vec J_{\vec j}}^{\ell+1})^{-1}
  -(\widetilde A_{h,\vec J_{\vec j}}^{\ell+1})^{-1}\right]
  \mathcal R_{h,\vec J_{\vec j}}^{\ell+1}
  A_{h,\vec J}^{\ell}v_{h,\vec J}\\
  &\quad=\sum_{\vec j\in\mathcal S}
  \widetilde{\mathcal R}_{h,\vec J_{\vec j}}^{\ell+1,\top}
  \left[I-(\widetilde A_{h,\vec J_{\vec j}}^{\ell+1})^{-1}
  A_{h,\vec J_{\vec j}}^{\ell+1}\right]
  (A_{h,\vec J_{\vec j}}^{\ell+1})^{-1}
  \mathcal R_{h,\vec J_{\vec j}}^{\ell+1}
  A_{h,\vec J}^{\ell}v_{h,\vec J}.
\end{aligned}
\]
{Lemma~\ref{lem:interpolation} and Assumption~\ref{ass:hkd}
yield a bound for the weighted prolongations.
Combining this bound with the local discrete inf--sup conditions,
continuity of the parent sesquilinear form, and the finite-overlap
estimate in~\cite[(2.10)]{graham2020domain} gives}
\[
\begin{aligned}
  &\left\|
  (\B_{h,\vec J}^{\ell}-\Bt_{h,\vec J}^{\ell})
  A_{h,\vec J}^{\ell}v_{h,\vec J}
  \right\|_{H^1_k(\Omega_{\vec J}^{\ell})}^2\\
  &\quad\lesssim\Lambda\sum_{\vec j\in\mathcal S}
  \left\|
  \left[I-(\widetilde A_{h,\vec J_{\vec j}}^{\ell+1})^{-1}
  A_{h,\vec J_{\vec j}}^{\ell+1}\right]
  (A_{h,\vec J_{\vec j}}^{\ell+1})^{-1}
  \mathcal R_{h,\vec J_{\vec j}}^{\ell+1}
  A_{h,\vec J}^{\ell}v_{h,\vec J}
  \right\|_{H^1_k(\Omega_{\vec J_{\vec j}}^{\ell+1})}^2\\
  &\quad\lesssim\Lambda\sum_{\vec j\in\mathcal S}
  \left\|
  I-(\widetilde A_{h,\vec J_{\vec j}}^{\ell+1})^{-1}
  A_{h,\vec J_{\vec j}}^{\ell+1}
  \right\|^2 (\mu_{h,\vec J_{\vec j}}^{\ell+1})^{-2}
  \|v_{h,\vec J}\|_{H^1_k(\Omega_{\vec J_{\vec j}}^{\ell+1})}^2\\
  &\quad\lesssim\Lambda^2
  \max_{\vec j\in\mathcal S}\left[
  (\mu_{h,\vec J_{\vec j}}^{\ell+1})^{-2}
  \left\|
  I-(\widetilde A_{h,\vec J_{\vec j}}^{\ell+1})^{-1}
  A_{h,\vec J_{\vec j}}^{\ell+1}
  \right\|^2\right]
  \|v_{h,\vec J}\|_{H^1_k(\Omega_{\vec J}^{\ell})}^2.
\end{aligned}
\]
Taking square roots proves the estimate.

\end{proof}

Next, we show that the local-solver errors can be controlled by choosing
the number of inner iterations $\mathcal N_\ell$ at each level appropriately, and then we show that the resulting nested preconditioner is a good approximation of the exact one-level operator.

\begin{lemma}[Levelwise block perturbation]
\label{lem:levelwise-block-perturbation}
Let $0\leq\ell<L-1$, and suppose that the levelwise estimates in
Lemma~\ref{lem:B-Bh} and Theorem~\ref{thm:bdd} hold uniformly. Then, for
every $\vec J\in\mathcal I_\ell$, we have
{
\[
\begin{aligned}
  &\left\|
  \left(I-\Bt_{h,\vec J}^{\ell}A_{h,\vec J}^{\ell}\right)^{\mathcal N}
  -\left(I-\B_{h,\vec J}^{\ell}A_{h,\vec J}^{\ell}\right)^{\mathcal N}
  \right\|\\
  &\qquad\lesssim
  \mathcal N
  \left(
    1+\Lambda
    \max_{\vec j\in\mathcal S}
    (\mu_{h,\vec J_{\vec j}}^{\ell+1})^{-1}
  \right)^{\mathcal N}
  (1+E_{\ell+1})^{\mathcal N-1}E_{\ell+1}.
\end{aligned}
\]
}
\end{lemma}
\begin{proof}
The power-difference identity 
$
  X^n-Y^n=\sum_{m=0}^{n-1}X^{n-1-m}(X-Y)Y^m
$ 
gives
\[
\begin{aligned}
  &\left\|
  \left(I-\Bt_{h,\vec J}^{\ell}A_{h,\vec J}^{\ell}\right)^{\mathcal N}
  -\left(I-\B_{h,\vec J}^{\ell}A_{h,\vec J}^{\ell}\right)^{\mathcal N}
  \right\|\\
  &\quad\leq
  \sum_{m=0}^{\mathcal N-1}
  \left\|
  I-\Bt_{h,\vec J}^{\ell}A_{h,\vec J}^{\ell}
  \right\|^{\mathcal N-1-m}
  \left\|
  (\B_{h,\vec J}^{\ell}-\Bt_{h,\vec J}^{\ell})
  A_{h,\vec J}^{\ell}
  \right\|
  \left\|
  I-\B_{h,\vec J}^{\ell}A_{h,\vec J}^{\ell}
  \right\|^m.
\end{aligned}
\]
Lemma~\ref{lem:B-Bh} and the definition of $E_{\ell+1}$ imply
\[
  \left\|
  (\B_{h,\vec J}^{\ell}-\Bt_{h,\vec J}^{\ell})
  A_{h,\vec J}^{\ell}
  \right\|
  \lesssim
  \Lambda
  \max_{\vec j\in\mathcal S}
  (\mu_{h,\vec J_{\vec j}}^{\ell+1})^{-1}E_{\ell+1}.
\]
The levelwise counterpart of Theorem~\ref{thm:bdd} gives
\[
  \left\|
  I-\B_{h,\vec J}^{\ell}A_{h,\vec J}^{\ell}
  \right\|
  \lesssim
  1+\Lambda
  \max_{\vec j\in\mathcal S}
  (\mu_{h,\vec J_{\vec j}}^{\ell+1})^{-1}.
\]
{
The triangle inequality and the preceding two estimates also give
\[
\begin{aligned}
  \left\|
  I-\Bt_{h,\vec J}^{\ell}A_{h,\vec J}^{\ell}
  \right\|
  &\lesssim
  \left(
    1+\Lambda
    \max_{\vec j\in\mathcal S}
    (\mu_{h,\vec J_{\vec j}}^{\ell+1})^{-1}
  \right)(1+E_{\ell+1}).
\end{aligned}
\]
}
Substituting these bounds into the power-difference estimate yields the result.
\end{proof}

\begin{lemma}[Error bounds for the nested local solves]
\label{lem:nested-local-errors}
Under the hypotheses of Theorem~\ref{thm:nested-conv}, the iteration
counts in~\eqref{eq:inner-iteration-count} yield the local solver error bounds
\[
  E_L=0,
  \qquad
  E_\ell
  \leq
  \left(\frac{1+\rho}{2}\right)^{\mathcal N_\ell/\mathcal N},
  \qquad 1\leq\ell<L.
\]
\end{lemma}
\begin{proof}
The terminal solves are exact, so $E_L=0$. If $L=1$, there are no nonterminal local levels, so the result reduces to
$E_1=0$. For $L\geq2$, we prove the remaining estimates by
downward induction.
At level $L-1$,
$\Bt_{h,\vec J}^{L-1}=\B_{h,\vec J}^{L-1}$ because the
level-$L$ child solves are exact. The definition of the approximate inverse $\tilde{A}_{h,\vec J}^{L-1}$ and
Assumption~\ref{ass:uniform-discrete-contraction} imply
\[
\begin{aligned}
  E_{L-1}
  =\max_{\vec J}\left\|
  I-(\widetilde A_{h,\vec J}^{L-1})^{-1}
  A_{h,\vec J}^{L-1}
  \right\|&=\max_{\vec J}\left\|
  \left(I-\B_{h,\vec J}^{L-1}A_{h,\vec J}^{L-1}\right)^{
  \mathcal N_{L-1}}\right\|\\
  &\leq
  \rho^{\mathcal N_{L-1}/\mathcal N}
  \leq
  \left(\frac{1+\rho}{2}\right)^{
  \mathcal N_{L-1}/\mathcal N}.
\end{aligned}
  \]
Let $1\leq\ell<L-1$ and assume that
$  E_{\ell+1}
  \leq
  \left(\frac{1+\rho}{2}\right)^{
  \mathcal N_{\ell+1}/\mathcal N}
  \leq1.$ 
Since $\mathcal N_\ell/\mathcal N$ is an integer
by~\eqref{eq:inner-iteration-count},
\eqref{eq:recursive-approximate-inverse} shows that it is enough to prove
contraction over one fixed block of $\mathcal N$ iterations, namely,
\begin{equation}\label{eq:nested-fixed-block-contraction}
  \begin{aligned}
    &~~\max_{\vec J}\left\|
  \left(I-\Bt_{h,\vec J}^{\ell}A_{h,\vec J}^{\ell}\right)^{\mathcal N}
  \right\|
  &\leq\frac{1+\rho}{2}.
  \end{aligned}
\end{equation}
From Assumption~\ref{ass:uniform-discrete-contraction}, we have
$
  \max_{\vec J}\left\|
  \left(I-\B_{h,\vec J}^{\ell}A_{h,\vec J}^{\ell}\right)^{\mathcal N}
  \right\|\leq\rho.
$ 
Thus, it remains to show that the second term in~\eqref{eq:nested-fixed-block-contraction} is bounded by $(1-\rho)/2$.
Since $E_{\ell+1}\leq1$ and $\mathcal N$ is fixed,
Lemma~\ref{lem:levelwise-block-perturbation}, the assumed local inf--sup
bound, and the uniform boundedness of $\Lambda$ give
\[
\begin{aligned}
  &\max_{\vec J}\left\|
  \left(I-\Bt_{h,\vec J}^{\ell}A_{h,\vec J}^{\ell}\right)^{\mathcal N}
  -\left(I-\B_{h,\vec J}^{\ell}A_{h,\vec J}^{\ell}\right)^{\mathcal N}
  \right\|\lesssim
  \left(1+\frac{k}{2^{\ell+1}}\right)^{\mathcal N}E_{\ell+1}.
\end{aligned}
\]
On the other hand, the induction hypothesis and
\eqref{eq:inner-iteration-count} imply
\[
\begin{aligned}
  E_{\ell+1}
  &\leq
  \left(\frac{1+\rho}{2}\right)^{
  \mathcal N_{\ell+1}/\mathcal N}\leq
  \left(\frac{1+\rho}{2}\right)^{q_0}
  \left(1+\frac{k}{2^{\ell+1}}\right)^{
  -q_1\log\!\left(2/(1+\rho)\right)}.
\end{aligned}
\]
Consequently,
\[
\begin{aligned}
  &\max_{\vec J}\left\|
  \left(I-\Bt_{h,\vec J}^{\ell}A_{h,\vec J}^{\ell}\right)^{\mathcal N}
  -\left(I-\B_{h,\vec J}^{\ell}A_{h,\vec J}^{\ell}\right)^{\mathcal N}
  \right\|\\
  &\qquad\lesssim
  \left(\frac{1+\rho}{2}\right)^{q_0}
  \left(1+\frac{k}{2^{\ell+1}}\right)^{
  \mathcal N-q_1\log\!\left(2/(1+\rho)\right)}.
\end{aligned}
\]
The choice
$
  q_1\log\!\left(\frac{2}{1+\rho}\right)>\mathcal N
$
makes the last exponent negative. The frequency-dependent factor is therefore
bounded by one, and $q_0$ can then be chosen sufficiently large, independently
of $k$, $L$, and $\ell$, so that
\[
  \max_{\vec J}\left\|
  \left(I-\Bt_{h,\vec J}^{\ell}A_{h,\vec J}^{\ell}\right)^{\mathcal N}
  -\left(I-\B_{h,\vec J}^{\ell}A_{h,\vec J}^{\ell}\right)^{\mathcal N}
  \right\|
  \leq\frac{1-\rho}{2}.
\]
Since $\mathcal N_\ell$ is an integer multiple of $\mathcal N$,
\eqref{eq:recursive-approximate-inverse}, the preceding estimate, and
Assumption~\ref{ass:uniform-discrete-contraction} yield
\[
\begin{aligned}
  E_\ell
  &\leq
  \left(
    \max_{\vec J}\left\|
    \left(I-\Bt_{h,\vec J}^{\ell}A_{h,\vec J}^{\ell}\right)^{\mathcal N}
    \right\|
  \right)^{\mathcal N_\ell/\mathcal N}\leq
  \left(\rho+\frac{1-\rho}{2}\right)^{
    \mathcal N_\ell/\mathcal N}
  =
  \left(\frac{1+\rho}{2}\right)^{
    \mathcal N_\ell/\mathcal N}.
\end{aligned}
\]
This closes the downward induction.
\end{proof}

\begin{proof}[{Proof of Theorem~\ref{thm:nested-conv}}]
If $L=1$, the result follows directly from
Assumption~\ref{ass:uniform-discrete-contraction}. {For $L\geq2$,
Lemmas~\ref{lem:nested-local-errors}
and~\ref{lem:levelwise-block-perturbation} with $\ell=0$, together with
the level-one inf--sup bound and the uniform boundedness of $\Lambda$, give
$E_1\leq1$ and
\[
  \left\|(I-\Bt_hA_h)^{\mathcal N}-(I-\B_hA_h)^{\mathcal N}\right\|
  \lesssim \left(1+\frac{k}{2}\right)^{\mathcal N}E_1.
\]
Moreover, \eqref{eq:inner-iteration-count} and
Lemma~\ref{lem:nested-local-errors} imply
\[
  E_1\leq\left(\frac{1+\rho}{2}\right)^{\mathcal N_1/\mathcal N}
  \leq\left(\frac{1+\rho}{2}\right)^{q_0}
  \left(1+\frac{k}{2}\right)^{-q_1\log(2/(1+\rho))}.
\]
Since $q_1\log(2/(1+\rho))>\mathcal N$, substituting this estimate for
$E_1$ into the preceding perturbation bound shows that $q_0$ can be chosen
sufficiently large, independently of $k$ and $L$, so that
\[
  \left\|(I-\Bt_hA_h)^{\mathcal N}-(I-\B_hA_h)^{\mathcal N}\right\|
  \leq\frac{1-\rho}{2}.
\]
}
Combining this with
Assumption~\ref{ass:uniform-discrete-contraction} gives
\[
  \left\|\left(I-\Bt_hA_h\right)^{\mathcal N}\right\|
  \leq\frac{1-\rho}{2}+\rho
  =\frac{1+\rho}{2}<1.
\]
\end{proof}

\section{Numerical experiments}
\label{sec:experiments}
This section assesses the convergence and parallel performance of the
RAS--PML method. {We first study how the inner iteration count $q$ affects
the convergence of Nested($q$), defined in
Section~\ref{ssec:nested-method}.
We then examine the choice of local PML and overlap widths
for one-level RAS--PML and compare RAS with PML and impedance
transmission conditions.} With the
selected parameters, we measure iteration counts and parallel runtimes
as the wavenumber and the number of subdomains increase together.
Finally, we examine sensitivity to mesh refinement and assess the method
on heterogeneous geological models.

Unless otherwise stated, we consider the constant-coefficient problem
with $c\equiv1$. The computational domain, including the global PML,
is the unit box $\Omega=(0,1)^d$, $d=2,3$.
We use the PML profile $g(x)=10kx^3$
and impose a homogeneous Dirichlet condition on the outer boundary
$\partial\Omega$. The forcing is a Gaussian centred at the centre
$\boldsymbol{x}_0$ of the computational domain,
\[
f(\boldsymbol{x}) =
C_d k^d
\exp\Bigl(-\frac{16 k^2|\boldsymbol{x}-\boldsymbol{x}_0|^2}{\pi^2}\Bigr),
\qquad C_2=\frac{16}{\pi^3},\quad C_3=\frac{16}{\pi^{9/2}}.
\]
{This Gaussian does not satisfy the support condition on $f$ imposed
in Section~2.1. It is strongly localized around $\boldsymbol{x}_0$,
with exponentially decaying tails.}

For the spatial discretization of the PML truncated problem, we use conforming
$Q_1$ finite elements on uniform Cartesian meshes with mesh size $h$.
We denote the number of complex global unknowns by DoF.
In the constant-coefficient frequency studies, 
{as motivated in the introduction, we study the case when there is a fixed number of grid points per wavelength.}
We denote this resolution by PPW, with wavelength
$\lambda = 2\pi/k$ and $\mathrm{PPW}=\lambda/h=2\pi/(kh)$.
For the domain decomposition, we use Cartesian partitions. An isotropic
partition with $n$ subdomains in each coordinate direction contains
$N=n^d$ subdomains in total. In these studies, we increase $n$ proportionally to $k$, {thus} keeping the number of
points in each nominal subdomain core fixed. The partition-of-unity
{functions $\{\chi_j\}$} are constructed as tensor products of one-dimensional
piecewise-linear functions on the overlaps.
The choice of the overlap width  $\delta$ and the local PML width 
$\kappa$ will be discussed in Section~\ref{sec:choice-pml}.

Unless otherwise stated, we use Richardson iteration with unit relaxation
parameter and a zero initial guess. We also use restarted GMRES with
restart length $30$, denoted by GMRES(30), in the experiments specified
below. For both methods, convergence is declared {when the relative residual in the Euclidean norm is below $10^{-10}$.} The maximum number of iterations is $2000$.
In the tables, $>2000$ indicates that this tolerance
was not reached within the iteration limit, while $\times$ denotes
divergence.

The reported times use mixed precision: local factorizations and solves
are performed by MUMPS in single precision, while global arithmetic is
performed in double precision. These single-precision local solves may
increase the iteration count by only a few steps. The computations are run on
nodes with two 48-core Intel Xeon Platinum 9242 processors and
$375$~GiB of memory. We assign one MPI rank to each subdomain
and use one thread per rank.
The reported times are labelled Setup for local factorization, {Solve time}, and Total for their sum.

\subsection{Convergence of the nested RAS--PML method}
\label{sec:nested-convergence}

We investigate the effect of the inner iteration count by applying
Algorithm~\ref{alg:nestedras} to~\eqref{eq:weak}, using a
recursive $2\times2$ decomposition with $L$ levels, $n=2^L$ subdomains
in each coordinate direction, and $N=(2^L)^2$ finest-level subdomains.

For the tests in Table~\ref{tab:nested_ras}, we set $kh=1/2$ and fix the
local PML and overlap widths at $\kappa=\pi/75$ and $\delta=1/75$,
independently of $k$. We compare Nested($q$) {for $q\in\{2,3\}$,}
where each non-terminal inner level uses $\mathcal N_\ell=q$ iterations.
The reported counts $a(b)$ give the number of outer iterations $a$
and the cumulative number of terminal-level iterations $b$ along a
single recursive path. These counts satisfy $b=aq^{L-1}$, since each
of the $L-1$ non-terminal inner levels contributes a factor of $q$.
{All three methods converge in every test. For $L>1$, Nested(2)
and Nested(3) require fewer outer iterations than one-level RAS--PML,
but more local solves per finest-level subdomain. 
We do not verify the equivalence in
Theorem~\ref{thm:nested-ras-equivalence} numerically.
For ease of implementation, the one-level experiments use standard  PoU weights rather than the
hierarchy-induced composite weights required by the theorem.
These results motivate the use of one-level RAS--PML in the remaining
experiments.}

\begin{table}[!ht]
    \centering
    \footnotesize
    \makeatletter
    \long\def\@makecaption#1#2{%
        \footnotesize\centering
        {\normalfont\scshape #1}. {\normalfont\itshape #2}\par
        \vskip\belowcaptionskip}
    \makeatother
    \caption{Iteration counts for nested RAS--PML with
    $kh=1/2$, 
    $\kappa=\pi/75$ and $\delta=1/75$. Entries $a(b)$ give outer
    iterations $a$ and cumulative terminal-level iterations $b$ along a single
    recursive path.}
    \label{tab:nested_ras}
    \setlength{\tabcolsep}{4pt}
    \renewcommand{\arraystretch}{1.05}
    \begin{tabular}{@{}cccccc@{}}
    \toprule
    \multicolumn{3}{c}{Problem setting}
    & \multicolumn{3}{c}{Iteration counts} \\
    \cmidrule(lr){1-3}\cmidrule(lr){4-6}
    $k$ & DoF ($10^6$) & $N$
    & Nested(3) & Nested(2) &One-level RAS-PML \\
    \midrule
    150  & 0.090 & $2^2$  & 8 (8)  & 8 (8)  & 8 \\
    300  & 0.360 & $4^2$  & 7 (21) & 9 (18) & 14 \\
    600  & 1.440 & $8^2$  & 5 (45) & 10 (40) & 24 \\
    1200 & 5.760 & $16^2$ & 4 (108) & 9 (72) & 36 \\
    \bottomrule
    \end{tabular}
\end{table}

\subsection{Choice of PML and overlap widths}
\label{sec:choice-pml}
\label{sec:3d}

We next examine the choice of local PML and overlap widths as the
wavenumber and {hence the} number of subdomains increase. We compare widths that
are fixed multiples of the wavelength with widths {which grow with $\log k$ times the wavelength as $k$ increases}. Following  \cite{xie_dd29_raspml}, we set $k_0=k/n$ and define $\ell(k):=\frac12\frac{k}{k-k_0}\log_2(k/k_0),~ k>k_0.$ 
The nominal local PML and overlap widths are {taken to be}
\begin{equation}
\label{eq:nominal-widths}
\kappa=c_\kappa\ell(k)\lambda,
\qquad
\delta=c_\delta\ell(k)\lambda,
\end{equation}
both rounded down to mesh layers.
We use $c_\delta=1/\pi$ in both dimensions, with $c_\kappa=1$ in two
dimensions and $c_\kappa=1/2$ in three.
At fixed $kh$, both {PML and overlap widths} contain $\mathcal O(\log k)$ mesh
layers and have physical width $\mathcal O(k^{-1}\log k)$.

{Tables~\ref{tab:bscc-parameter-ablation} and~\ref{tab:3d-parameter-choice}
report Richardson iteration counts in two dimensions with $kh=1/2$ and
in three dimensions with $kh=1$, respectively. We compare three choices:
both widths fixed in wavelength units; only the PML width scaled by
$\ell(k)$; and both widths scaled by $\ell(k)$. In two dimensions, fixed
widths give 63 iterations at $k=1200$, followed by Richardson divergence
at all tested $k\ge2400$. Scaling the PML width restores convergence
through $k=9600$. Scaling the overlap as well reduces the count from
305 to 232 at this wavenumber. The three-dimensional results show a
similar trend. Over the tested range, increasing the PML width in
wavelength units maintains convergence. Increasing the overlap as well
further reduces iteration counts.}

\begin{table}[!ht]
\centering
\footnotesize
\makeatletter
\long\def\@makecaption#1#2{%
    \footnotesize\centering
    {\normalfont\scshape #1}. {\normalfont\itshape #2}\par
    \vskip\belowcaptionskip}
\makeatother
\caption{Richardson iteration counts in 2D with $kh=1/2$, $k_0=150$ and various PML and overlap widths.}
\label{tab:bscc-parameter-ablation}
\setlength{\tabcolsep}{4pt}
\renewcommand{\arraystretch}{1.05}
\begin{tabular}{@{}cccccc@{}}
\toprule
\multicolumn{3}{c}{Problem setting} & \multicolumn{3}{c}{Richardson iteration counts for $(\kappa,\delta)$} \\
\cmidrule(lr){1-3}\cmidrule(lr){4-6}
$k$ & DoF ($10^6$) & $N$
& $(\lambda,\lambda/\pi)$
& $(\ell(k)\lambda,\lambda/\pi)$
& $(\ell(k)\lambda,\ell(k)\lambda/\pi)$ \\
\midrule
300 & 0.360 & $2^2$ & 7 & 7 & 7 \\
600 & 1.440 & $4^2$ & 16 & 15 & 14 \\
1200 & 5.760 & $8^2$ & 63 & 33 & 32 \\
2400 & 23.040 & $16^2$ & $\times$ & 74 & 66 \\
4800 & 92.160 & $32^2$ & $\times$ & 154 & 125 \\
9600 & 368.640 & $64^2$ & $\times$ & 305 & 232 \\
\bottomrule
\end{tabular}
\end{table}

\begin{table}[!ht]
\centering
\footnotesize
\makeatletter
\long\def\@makecaption#1#2{%
    \footnotesize\centering
    {\normalfont\scshape #1}. {\normalfont\itshape #2}\par
    \vskip\belowcaptionskip}
\makeatother
\caption{Richardson iteration counts in 3D with $kh=1$, $k_0=30$ and various PML and overlap widths.}
\label{tab:3d-parameter-choice}
\setlength{\tabcolsep}{4pt}
\renewcommand{\arraystretch}{1.05}
\begin{tabular}{@{}cccccc@{}}
\toprule
\multicolumn{3}{c}{Problem setting} & \multicolumn{3}{c}{Richardson iteration counts for $(\kappa,\delta)$} \\
\cmidrule(lr){1-3}\cmidrule(lr){4-6}
$k$ & DoF ($10^6$) & $N$
& $(\lambda/2,\lambda/\pi)$
& $(\ell(k)\lambda/2,\lambda/\pi)$
& $(\ell(k)\lambda/2,\ell(k)\lambda/\pi)$ \\
\midrule
120 & 1.728 & $4^3$ & 32 & 23 & 23 \\
240 & 13.824 & $8^3$ & $>2000$ & 46 & 44 \\
480 & 110.592 & $16^3$ & $\times$ & 97 & 74 \\
600 & 216.000 & $20^3$ & $\times$ & 125 & 98 \\
\bottomrule
\end{tabular}
\end{table}

\subsection{Comparison of transmission conditions}
\label{sec:transmission-comparison}

Using the logarithmic PML {and overlap} widths in the last columns of
Tables~\ref{tab:bscc-parameter-ablation} and~\ref{tab:3d-parameter-choice} in two and three dimensions, respectively, we compare RAS--PML with
RAS--Imp, which uses the impedance
condition on artificial boundaries~\cite{gong2022convergence}.
Both methods use the same global problem, non-overlapping
partition and overlap width $\ell(k)\lambda/\pi$.

Table~\ref{tab:krylov-comparison} compares RAS--PML and RAS--Imp
using GMRES(30). Both preconditioners yield convergence at every tested
wavenumber, but RAS--PML requires fewer iterations and shorter solve
times in every case, at the expense of longer local factorization
times. Thus, PML transmission improves convergence and reduces the
iterative solution time relative to impedance transmission in these tests.
For RAS--PML, Table~\ref{tab:krylov-comparison} also compares GMRES(30)
with Richardson. GMRES(30) does not consistently reduce the iteration
count: in two dimensions, the counts
are nearly unchanged at the two lower wavenumbers and increase at
$k=9600$, while solve times {for GMRES(30)} are longer in all three cases. In three
dimensions, GMRES(30) reduces the counts at the two higher wavenumbers,
but the time savings remain modest.
Using GMRES without restarting would require storing a Krylov basis
that grows with the iteration count~\cite{saad1986gmres}. Given this memory cost and the
limited gains from GMRES(30), we use Richardson for the
constant-coefficient frequency study below.

\begin{table}[!ht]
\centering
\footnotesize
\makeatletter
\long\def\@makecaption#1#2{%
    \footnotesize\centering
    {\normalfont\scshape #1}. {\normalfont\itshape #2}\par
    \vskip\belowcaptionskip}
\makeatother
\caption{Iteration counts and wall-clock times (seconds) for RAS--Imp and RAS--PML with $kh=1/2$, $k_0=150$ in 2D, and $kh=1$, $k_0=30$ in 3D}
\label{tab:krylov-comparison}
\setlength{\tabcolsep}{4pt}
\renewcommand{\arraystretch}{1.05}
\begin{tabular}{@{}cccccccccccc@{}}
\toprule
& & & \multicolumn{3}{c}{RAS--Imp} & \multicolumn{6}{c}{RAS--PML} \\
\cmidrule(lr){4-6}\cmidrule(lr){7-12}
& & & \multicolumn{3}{c}{GMRES(30)} & \multicolumn{3}{c}{GMRES(30)} & \multicolumn{3}{c}{Richardson} \\
\cmidrule(lr){4-6}\cmidrule(lr){7-9}\cmidrule(lr){10-12}
Dimension & $k$ & $N$ & Iter. & Setup & Solve & Iter. & Setup & Solve & Iter. & Setup & Solve \\
\midrule
\multirow{3}{*}{2D} & 1200 & $8^2$ & 96 & 1.13 & 9.76 & 30 & 1.28 & 3.53 & 32 & 1.30 & 3.12 \\
 & 4800 & $32^2$ & 430 & 1.38 & 64.78 & 124 & 1.89 & 23.02 & 125 & 1.86 & 18.66 \\
 & 9600 & $64^2$ & 958 & 1.52 & 185.34 & 322 & 2.07 & 78.97 & 232 & 2.10 & 41.66 \\
\midrule
\multirow{3}{*}{3D} & 120 & $4^3$ & 55 & 2.18 & 3.48 & 23 & 3.70 & 2.20 & 23 & 3.70 & 2.07 \\
 & 300 & $10^3$ & 129 & 4.34 & 16.78 & 48 & 8.25 & 11.52 & 52 & 8.26 & 12.17 \\
 & 600 & $20^3$ & 296 & 6.98 & 74.55 & 77 & 15.96 & 38.92 & 98 & 15.95 & 39.71 \\
\bottomrule
\end{tabular}
\end{table}

\subsection{\texorpdfstring{$\mathcal{O}(k)$ performance with growing frequency}{O(k) performance with growing frequency}}
\label{sec:large-scale}
\label{sec:time-analysis}

We now examine how the iteration counts and parallel runtimes grow
with $k$ for {the logarithmic PML and overlap} widths.
We fix $kh$, using
$(kh,k_0)=(1/2,150)$ in two dimensions and $(1,30)$ in three.
The number of points in each nominal core then remains fixed,
while the global number of unknowns and subdomains grows as $k^d$.

Tables~\ref{tab:bscc-scaling-2d} and~\ref{tab:bscc-scaling-3d}
report the Richardson iteration counts and parallel timings, with
convergence at every tested wavenumber in both dimensions.
The largest two-dimensional system has $1.47\times10^9$ complex
unknowns on $128^2$ subdomains; the largest three-dimensional system
has $4.75\times10^8$ unknowns on $26^3$ subdomains.

\begin{table}[!ht]
\centering
\footnotesize
\makeatletter
\long\def\@makecaption#1#2{%
    \footnotesize\centering
    {\normalfont\scshape #1}. {\normalfont\itshape #2}\par
    \vskip\belowcaptionskip}
\makeatother
\caption{Richardson iteration counts and wall-clock times in two dimensions
for $kh=1/2$, $\kappa=\ell(k)\lambda$,
$\delta=\ell(k)\lambda/\pi$ and $k_0=150$.
For $k=4800$, $9600$, and $19200$, the global residual norm is evaluated
only at termination.}
\label{tab:bscc-scaling-2d}
\setlength{\tabcolsep}{4pt}
\renewcommand{\arraystretch}{1.05}
\begin{tabular}{@{}cccccccc@{}}
\toprule
\multicolumn{3}{c}{Problem setting} &
\multicolumn{5}{c}{Richardson iteration} \\
\cmidrule(lr){1-3}\cmidrule(lr){4-8}
$k$ & DoF ($10^6$) & $N$ & Iter. & Setup (s) &
Solve (s) & Total (s) & Total/$k$ ($10^{-3}$) \\
\midrule
300   & 0.360    & 4     & 7   & 0.85 & 0.39  & 1.24  & 4.13 \\
600   & 1.440    & 16    & 14  & 1.02 & 0.93  & 1.96  & 3.27 \\
1200  & 5.760    & 64    & 32  & 1.34 & 3.16  & 4.50  & 3.75 \\
2400  & 23.040   & 256   & 66  & 1.60 & 8.21  & 9.88  & 4.12 \\
4800  & 92.160   & 1024  & 125 & 1.88 & 16.62 & 20.50 & 4.27 \\
9600  & 368.640  & 4096  & 232 & 2.10 & 33.65 & 35.75 & 3.72 \\
19200 & 1474.560 & 16384 & 485 & 2.58 & 79.88 & 82.46 & 4.30 \\
\bottomrule
\end{tabular}
\end{table}

\begin{table}[!ht]
\centering
\footnotesize
\makeatletter
\long\def\@makecaption#1#2{%
    \footnotesize\centering
    {\normalfont\scshape #1}. {\normalfont\itshape #2}\par
    \vskip\belowcaptionskip}
\makeatother
\caption{Richardson iteration counts and wall-clock times in 3D with
$kh=1$, $\kappa=\ell(k)\lambda/2$, $\delta=\ell(k)\lambda/\pi$ and $k_0=30$.
At $k=480$ and $k=780$, the global residual norm is evaluated only
at termination.}
\label{tab:bscc-scaling-3d}
\setlength{\tabcolsep}{4pt}
\renewcommand{\arraystretch}{1.05}
\begin{tabular}{@{}cccccccc@{}}
\toprule
\multicolumn{3}{c}{Problem setting} &
\multicolumn{5}{c}{Richardson iteration} \\
\cmidrule(lr){1-3}\cmidrule(lr){4-8}
$k$ & DoF ($10^6$) & $N$ & Iter. & Setup (s) &
Solve (s) & Total (s) & Total/$k$ ($10^{-3}$) \\
\midrule
60  & 0.216   & 8     & 10  & 1.65  & 0.42  & 2.07  & 34.50 \\
120 & 1.728   & 64    & 23  & 3.70  & 2.08  & 5.76  & 48.00 \\
240 & 13.824  & 512   & 44  & 7.59  & 9.19  & 16.95 & 70.63 \\
480 & 110.592 & 4096  & 74  & 12.94 & 22.13 & 35.07 & 73.06 \\
780 & 474.552 & 17576 & 110 & 17.15 & 40.21 & 57.36 & 73.54 \\
\bottomrule
\end{tabular}
\end{table}

\begin{figure}[!htbp]
    \centering
    \begin{minipage}{0.4\textwidth}
        \centering
        \includegraphics[width=\linewidth]{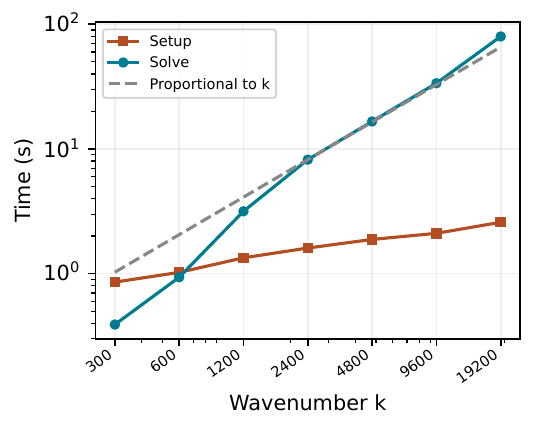}
        \par\smallskip{\small (e) Two-dimensional setup and solve times.}\par
    \end{minipage}\hspace{0.04\textwidth}
    \begin{minipage}{0.4\textwidth}
        \centering
        \includegraphics[width=\linewidth]{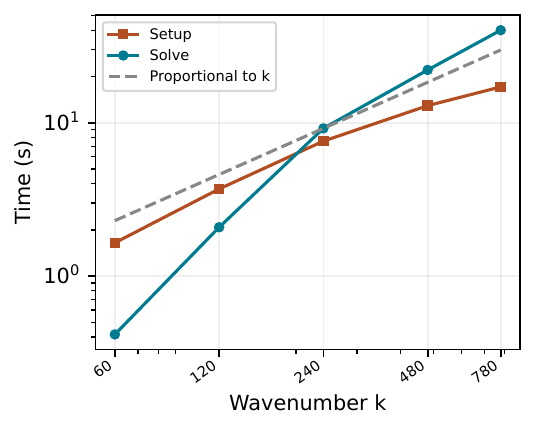}
        \par\smallskip{\small (f) Three-dimensional setup and solve times.}\par
    \end{minipage}
    \caption{Frequency dependence of Richardson convergence and parallel
    timings. Panels (a) and (b) show iteration counts divided by
    $n=k/k_0$. Panels (c) and (d) show relative residual histories
    against $m/n$. Panels (e) and (f) show Setup and Solve times;
    dashed lines indicate growth proportional to $k$.}
    \label{fig:bscc-frequency-scaling}
    \label{fig:normalized-frequency}
\end{figure}

The iteration counts divided by $n=k/k_0$ in
Figure~\ref{fig:normalized-frequency} lie between $3.50$ and $4.13$
in two dimensions and between $4.23$ and $5.75$ in three.
The residual histories also
approximately coincide when plotted against $m/n$, where $m$ is the
iteration index, particularly in two dimensions. Both observations
are consistent with iteration counts proportional to $n$, and hence
to $k$ in these experiments, while the global number of unknowns grows
as $k^d$. The solve times also grow approximately linearly
with $k$ at the higher wavenumbers, whereas factorization time grows
more slowly. {For $k\ge4800$ in two dimensions and $k\ge480$ in three,
the timed runs use iteration counts determined in separate runs that
check the relative residual at every iteration. In the timed runs,
the relative residual is evaluated only at termination and is below
$10^{-10}$ in every case.} Together, the convergence and timing results support
approximately $\mathcal O(k)$ performance over the tested range.

\begingroup
\capstartfalse 
\begin{figure}[!htbp]
    \centering
    \begin{minipage}{0.4\textwidth}
        \centering
        \includegraphics[width=\linewidth]{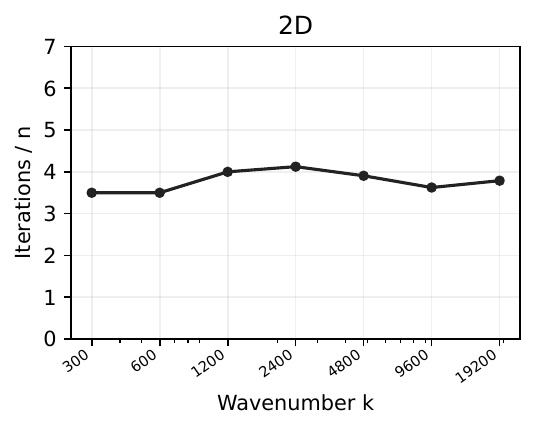}
        \par\smallskip{\small (a) Two-dimensional iteration counts.}\par
    \end{minipage}\hspace{0.04\textwidth}
    \begin{minipage}{0.4\textwidth}
        \centering
        \includegraphics[width=\linewidth]{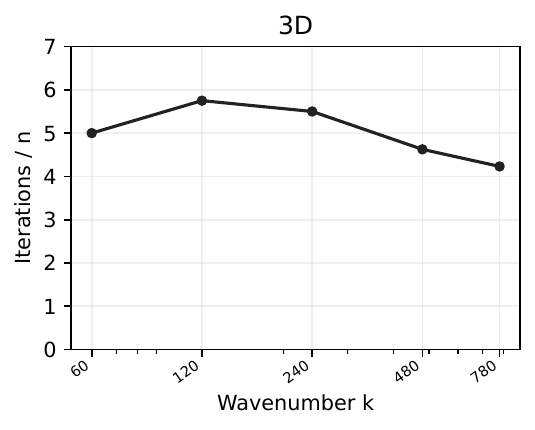}
        \par\smallskip{\small (b) Three-dimensional iteration counts.}\par
    \end{minipage}
\end{figure}

\begin{figure}[!htbp]
    \centering
    \begin{minipage}{0.4\textwidth}
        \centering
        \includegraphics[width=\linewidth]{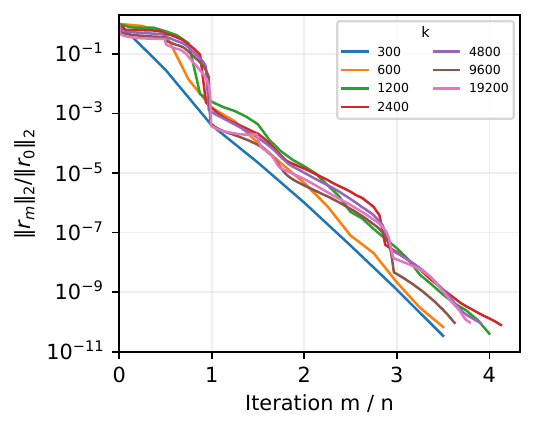}
        \par\smallskip{\small (c) Two-dimensional residual histories.}\par
    \end{minipage}\hspace{0.04\textwidth}
    \begin{minipage}{0.4\textwidth}
        \centering
        \includegraphics[width=\linewidth]{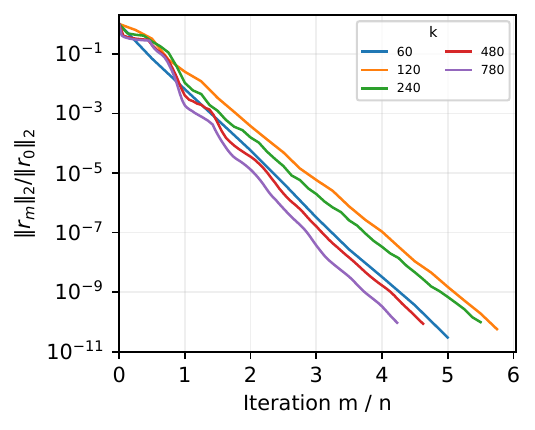}
        \par\smallskip{\small (d) Three-dimensional residual histories.}\par
    \end{minipage}
\end{figure}
\endgroup

\FloatBarrier

\subsection{Convergence under mesh refinement}
\label{sec:mesh-resolution}


Table~\ref{tab:ras-pml-mesh} compares the coarse mesh $kh=1$, the
standard mesh $kh=1/2$, and {finer meshes} with
$h^{-1}\approx0.12k^{3/2}$.
{For nontrapping Helmholtz problems with sufficient regularity,
analyses with radial PML truncation provide first-order FEM error
estimates when $k^3h^2$ is sufficiently small~\cite{li2019fem,Galkowski_2025}.
This motivates the $h\sim k^{-3/2}$ scaling.
}
{The iteration counts are nearly unchanged across the tested meshes.}

\begin{table}[!ht]
\centering
\footnotesize
\makeatletter
\long\def\@makecaption#1#2{%
    \footnotesize\centering
    {\normalfont\scshape #1}. {\normalfont\itshape #2}\par
    \vskip\belowcaptionskip}
\makeatother
\caption{Richardson iteration counts on coarse ($kh=1$), standard ($kh=1/2$)
and finer ($h^{-1}\approx0.12k^{3/2}$) meshes, with $k_0=150$
and the same physical PML and overlap widths
within each row. A dash denotes a configuration not tested.}
\label{tab:ras-pml-mesh}
\setlength{\tabcolsep}{4pt}
\renewcommand{\arraystretch}{1.05}
\begin{tabular}{@{}ccccc@{}}
\toprule
\multicolumn{2}{c}{Problem setting} & \multicolumn{3}{c}{Richardson iteration counts} \\
\cmidrule(lr){1-2}\cmidrule(lr){3-5}
$k$ & $N$ & $kh=1$ & $kh=1/2$ & $h^{-1}\approx0.12k^{3/2}$ \\
\midrule
300 & 4 & 8 & 7 & 7 \\
600 & 16 & 16 & 14 & 14 \\
1200 & 64 & 32 & 32 & 32 \\
2400 & 256 & 65 & 62 & 63 \\
4800 & 1024 & 126 & 123 & 122 \\
9600 & 4096 & 230 & 229 & -- \\
\bottomrule
\end{tabular}
\end{table}

\FloatBarrier

\subsection{Heterogeneous geological applications}
\label{sec:application}

We next consider media with strong velocity contrasts, faults and folded
layers. The Marmousi experiments examine whether the approximately
linear growth of Richardson iteration counts persists in a
two-dimensional heterogeneous medium. The three-dimensional
GO\_3D\_OBS and Overthrust experiments compare PML and impedance
transmission under GMRES(30).
{The model data specify the compressional-wave speed,
also called the P-wave speed, which we use as the wave speed in the
scalar acoustic equation~\eqref{eq:helmholtz0}.}

Marmousi contains folded layers and faults over a
$9.2\times3.0$~km region, with velocities from $1.5$ to $5.5$~km/s
\cite{versteeg1994marmousi,harikumar2025marmousi}.
GO\_3D\_OBS represents a crustal-scale subduction zone inspired by the
Nankai Trough, covering $175\times100\times30$~km with velocities from
$1.5$ to $8.634$~km/s~\cite{gorszczyk2021go,gorszczyk2020godata}.
The SEG/EAGE Overthrust model contains folded and thrust-faulted
structures over $20\times20\times5.1$~km, with
velocities from $1.5$ to $6.0$~km/s
\cite{aminzadeh1997models,aminzadeh1997overthrust}.
{
Let $v(\boldsymbol X)$ denote this physical wave-speed field at the
physical position $\boldsymbol X$. Figure~\ref{fig:bscc-3d-fields}(a,c)
shows sections of $v(\boldsymbol X)$ for GO\_3D\_OBS and Overthrust;
panels (b,d) show the real parts of the computed wavefields.
We introduce dimensionless coordinates $\boldsymbol x=\boldsymbol X/L_x$,
where $L_x$ is the horizontal model length. The dimensionless wave
speed $c(\boldsymbol x)$ in~\eqref{eq:helmholtz0} and the physical
frequency ${\nu}$ are then given by
\[
 c(\boldsymbol x)=\frac{v(L_x\boldsymbol x)}{v_{\rm ref}},\qquad
 {\nu}=\frac{k v_{\rm ref}}{2\pi L_x}.
\]
}
Here $v_{\rm ref}=1500$~m/s and ${\nu}$ is the physical frequency in Hz.
The local wavenumber is $k/c(\boldsymbol x)$. At the minimum velocity,
the Marmousi mesh provides $12.6$ PPW, and the
three-dimensional comparisons use $kh=1$, hence $6.3$ PPW.
The total number of subdomains grows approximately as
{$k^d$, where $d$ is the spatial dimension}.

The local PML and overlap {are given by}~\eqref{eq:nominal-widths}, evaluated at
$n=n_x=k/k_0$ in every coordinate direction, so
$\kappa=\ell(k)\lambda/2$ and $\delta=\ell(k)\lambda/\pi$.
In three dimensions the partition places about $30$ grid points along each
direction, so $n_x=k/30$ and $n_y,n_z$ follow the mesh.
GO\_3D\_OBS uses the grids $150\times86\times26$, $300\times171\times51$,
$600\times343\times103$ and $1190\times680\times204$, partitioned into
$5\times3\times1$, $10\times6\times2$, $20\times11\times3$ and
$40\times23\times7$ subdomains.
Overthrust uses $120\times120\times31$, $240\times240\times61$,
$480\times480\times122$ and $960\times960\times245$, partitioned into
$4\times4\times1$, $8\times8\times2$, $16\times16\times4$ and
$32\times32\times8$.
The global PML has nominal
thickness $3\lambda$, $1.5\lambda$ and $2\lambda$ for Marmousi,
GO\_3D\_OBS and Overthrust, respectively.

Richardson converges for Marmousi {up to} $k=9600$, with
iteration counts {divided by $n_x$ lying between $6.36$ and $9.25$.}
The approximately proportional iteration growth observed for the
constant-coefficient problem therefore persists over this frequency
range (Table~\ref{tab:marmousi-sweep}). At $k=9600$, the method solves
$1.20\times10^8$ complex unknowns on $64\times24$ subdomains in $407$
iterations and $61.84$ s. Factorization time remains small relative
to solve time, which therefore dominates the total. {Thus, up to $k = 9600$, the total solve time is roughly proportional to $k$.}
At $k=19200$, Richardson diverges with the same PML and overlap widths, whereas
GMRES(30) with the same RAS--PML preconditioner converges in $800$
steps, as marked by the asterisk in Table~\ref{tab:marmousi-sweep}.

\begin{table}[!ht]
\centering
\footnotesize
\makeatletter
\long\def\@makecaption#1#2{%
    \footnotesize\centering
    {\normalfont\scshape #1}. {\normalfont\itshape #2}\par
    \vskip\belowcaptionskip}
\makeatother
\caption{RAS--PML iteration counts and parallel timings for Marmousi, with $kh=1/2$ ($12.6$ PPW at the minimum velocity). The partitions satisfy $n_x=k/150$, $n_y=\lceil3n_x/8\rceil$ and $N=n_xn_y$. Nominal widths use $\kappa=\ell(k)\lambda$, $\delta=\ell(k)\lambda/\pi$ and $k_0=150$. Richardson is used except for the starred GMRES(30) entry.}
\label{tab:marmousi-sweep}
\setlength{\tabcolsep}{4pt}
\renewcommand{\arraystretch}{1.05}
\begin{tabular}{@{}ccccccccc@{}}
\toprule
\multicolumn{4}{c}{Problem setting} & \multicolumn{5}{c}{RAS--PML} \\
\cmidrule(lr){1-4}\cmidrule(lr){5-9}
$k$ & $\nu$ (Hz) & DoF ($10^6$) & $N$ & Iter. & Iter./$n_x$ & Setup (s) & Solve (s) & Total (s) \\
\midrule
600 & 15.57 & 0.469 & 8 & 30 & 7.50 & 0.60 & 1.25 & 1.85 \\
1200 & 31.14 & 1.877 & 24 & 74 & 9.25 & 1.09 & 4.88 & 5.98 \\
1800 & 46.71 & 4.223 & 60 & 106 & 8.83 & 1.09 & 8.60 & 9.70 \\
2400 & 62.28 & 7.512 & 96 & 134 & 8.38 & 1.41 & 14.13 & 15.54 \\
4800 & 124.56 & 30.048 & 384 & 259 & 8.09 & 1.51 & 31.13 & 32.63 \\
9600 & 249.11 & 120.192 & 1536 & 407 & 6.36 & 1.77 & 60.08 & 61.84 \\
19200 & 498.22 & 480.806 & 6144 & $800^{*}$ & $6.25^{*}$ & 2.24 & 223.19 & 225.43 \\
\bottomrule
\end{tabular}
\end{table}

For the three-dimensional models, we compare the preconditioners using
GMRES(30). Both methods use the same global PML, partition and {local PML/overlap widths}.
RAS--Imp adds no local PML and imposes the impedance condition with the
local wavenumber $k/c(\boldsymbol x)$ on artificial boundaries.

\begin{table}[htbp]
\centering
\footnotesize
\makeatletter
\long\def\@makecaption#1#2{%
    \footnotesize\centering
    {\normalfont\scshape #1}. {\normalfont\itshape #2}\par
    \vskip\belowcaptionskip}
\makeatother
\caption{GMRES(30) comparisons on three-dimensional geological models. Loc.\ is the largest local system, including overlap and local PML.}
\label{tab:bscc-applications}
\setlength{\tabcolsep}{2.4pt}
\renewcommand{\arraystretch}{1.05}
{\scriptsize
\begin{tabular}{@{}rrrrrrrrrrrr@{}}
\toprule
\multicolumn{4}{c}{Problem} & \multicolumn{4}{c}{RAS--PML} & \multicolumn{4}{c}{RAS--Imp} \\
\cmidrule(lr){1-4}\cmidrule(lr){5-8}\cmidrule(lr){9-12}
$k$ & $\nu$ (Hz) & DoF ($10^6$) & $N$ & Loc.\ ($10^3$) & Iter. & Setup (s) & Solve (s)
& Loc.\ ($10^3$) & Iter. & Setup (s) & Solve (s) \\
\midrule
\multicolumn{12}{@{}l}{GO\_3D\_OBS} \\
\addlinespace[2pt]
150 & 0.20 & 0.335 & 15 & 38.5 & 40 & 1.50 & 1.66 & 30.9 & 34 & 1.05 & 1.14 \\
300 & 0.41 & 2.616 & 120 & 64.3 & 56 & 4.59 & 7.71 & 42.9 & 59 & 2.24 & 4.60 \\
600 & 0.82 & 21.197 & 660 & 129.9 & 70 & 15.00 & 25.22 & 83.2 & 129 & 6.50 & 21.36 \\
1200 & 1.64 & 165.077 & 6440 & 160.2 & 242 & 33.10 & 92.83 & 99.3 & 458 & 8.38 & 85.02 \\
\midrule
\multicolumn{12}{@{}l}{Overthrust} \\
\addlinespace[2pt]
120 & 1.43 & 0.446 & 16 & 47.2 & 28 & 2.12 & 1.54 & 38.0 & 32 & 1.43 & 1.42 \\
240 & 2.86 & 3.514 & 128 & 75.5 & 59 & 5.48 & 7.65 & 51.7 & 105 & 2.86 & 8.15 \\
480 & 5.73 & 28.109 & 1024 & 110.4 & 90 & 13.49 & 27.07 & 74.0 & 197 & 5.49 & 32.68 \\
960 & 11.46 & 225.792 & 8192 & 169.3 & 168 & 20.22 & 58.18 & 106.0 & 425 & 8.95 & 77.48 \\
\bottomrule
\end{tabular}}
\end{table}

Both methods reach the tolerance at every frequency in
Table~\ref{tab:bscc-applications}. The iteration counts are close at the
lowest frequency of each model. At the highest frequency, RAS--PML takes
242 iterations on GO\_3D\_OBS and 168 on Overthrust, compared with 458
and 425 for RAS--Imp. The local PML makes the RAS--PML factorization
more expensive. On GO\_3D\_OBS the largest local system grows from
$38.5\times10^3$ to $160.2\times10^3$ unknowns for RAS--PML and from
$30.9\times10^3$ to $99.3\times10^3$ for RAS--Imp, while the global
problem grows from $0.335\times10^6$ to $165.077\times10^6$ unknowns.
These local systems remain
large relative to $\mathrm{DoF}/N$ because the model is thin in $z$ and
the overlap and local PML are taken from $n_x$. The sum of the
reported setup and solve times is shorter
for RAS--Imp except on Overthrust at $k=960$.

\begin{figure}[!htbp]
    \centering
    \capstart
    \input{figures/fig_bscc-3d-fields.tex}
\end{figure}


\section{Conclusion}
\label{sec:conclusion}
We have developed a one-level RAS--PML solver for high-frequency Helmholtz
problems in two and three dimensions, with approximately $\mathcal{O}(k)$
iteration counts and wall-clock time observed in constant-coefficient
experiments with $\mathcal{O}(k^d)$ concurrent local solves for systems
with $\mathcal{O}(k^d)$ unknowns.
{We prove that Nested(1) is algebraically equivalent to one-level
RAS--PML when the hierarchy-induced local problems and transfer
operators are retained. This equivalence provides a hierarchical
interpretation of the resulting one-level method.}
Moreover, we prove the convergence of the nested RAS--PML iteration under
suitable stability and contraction conditions with sufficiently accurate
recursive local solves.
Numerical experiments demonstrate
the effectiveness of the proposed parameter choices, the efficiency
of the algorithm, and its applicability to heterogeneous geological models.
In future work, we will build on the present $Q_1$ experiments to explore higher-order
discretizations and address the additional challenges they introduce,
motivated by the $hp$ framework of Graham and Spence~\cite{graham2025two}.

\section*{Acknowledgments}
We are grateful to Euan A. Spence (University of Bath) for valuable
comments and suggestions that improved this paper. We also thank
Pierre-Henri Tournier (Sorbonne University, Paris) and Scott Congreve
(Charles University, Prague) for very helpful discussions on hierarchical
one-level methods.

The authors used ChatGPT to improve the language and readability of the manuscript. The authors reviewed and edited all AI-assisted text and assume responsibility for all content.

SG was supported by the National Natural Science Foundation of China
(Grant Nos.~12671480 and 12411530066).

\FloatBarrier

\bibliographystyle{siamplain}
\bibliography{references}
\end{document}